\documentclass[12pt,draftclsnofoot,journal,letterpaper,onecolumn]{IEEEtran}

\usepackage{graphicx, wrapfig}
\usepackage{booktabs}
\usepackage{multirow}
\usepackage{dsfont}

\usepackage{hyperref}
\usepackage{color}
\usepackage{cite}
\usepackage{amsfonts,amsmath,amssymb,amsthm}
\usepackage[utf8]{inputenc}

\usepackage{float}
\usepackage{amsfonts,amsmath}
\usepackage{color}
\usepackage{multirow}
\usepackage{url}
\usepackage{array}
\usepackage{bbm}
\usepackage{bm}
\usepackage{enumitem}

\usepackage{epsfig}
\usepackage{graphicx, wrapfig}
\usepackage{booktabs}
\usepackage{multirow}

\usepackage{booktabs}
\usepackage[ruled,vlined]{algorithm2e}

\def\cE{{\cal E}}
\def\cM{{\cal M}}
\def\cC{{\cal C}}

\def\bgamma{\bar{\gamma}}
\def\cT{{\cal T}}

\def\cX{{\cal X}}

\def\lf{{\hbox{\tiny\rm lf}}}
\def\rg{{\hbox{\tiny\rm rg}}}
\def\tr{{\hbox{\rm\tiny tr}}}
\def\ts{{\hbox{\rm\tiny ts}}}

\def\Diag{{\hbox{\rm Diag}}}
 
\def\bR{{\mathbb{R}}}

\def\bE{{\mathbb{E}}}
\def\bP{{\mathbb {P}}}

\def\qed{\hfill$\Box$}

\def\Tr{{\hbox{\rm Tr}}}

\newcommand{\be}{\begin{eqnarray}}
\newcommand{\ee}[1]{\label{#1}\end{eqnarray}}

\newcommand{\ese}{\end{eqnarray*}}
\newcommand{\bse}{\begin{eqnarray*}}

\definecolor{MyDarkBlue}{rgb}{0,0.08,0.45}
\definecolor{MyViolet}{rgb}{0.45,0.08,0.95}
\definecolor{MyBrown}{rgb}{0.45,0.08,0}

\def\md{{\hbox{\tiny\rm md}}}
\def\ind{{\mathbbm{1}}}  % indicator function

\newcommand{\hide}[1]{}
\def\bR{{\mathbb{R}}}

\newtheorem{lemma}{\textbf{Lemma}}
\newtheorem{theorem}{\textbf{Theorem}}
\newtheorem{proposition}{\textbf{Proposition}}
\newtheorem{remark}{\textbf{Remark}}

\title{Sequential Change Detection by Optimal Weighted $\ell_2$ Divergence}
\author{Liyan~Xie and Yao~Xie
\thanks{Liyan~Xie and Yao~Xie are with School of Industrial and Systems Engineering, Georgia Institute of Technology,
Atlanta, GA, 30332 USA e-mail: {lxie49@gatech.edu, yao.xie@isye.gatech.edu}.}%
}

\begin{document}

\maketitle

\begin{abstract}
We present a new non-parametric statistic, called the weighed $\ell_2$ divergence, based on empirical distributions for sequential change detection. We start by constructing the weighed $\ell_2$ divergence as a fundamental building block for two-sample tests and change detection. The proposed statistic is proved to attain the optimal sample complexity in the offline setting. We then study the sequential change detection using the weighed $\ell_2$ divergence and characterize the fundamental performance metrics, including the average run length (ARL) and the expected detection delay (EDD). We also present practical algorithms to find the optimal projection to handle high-dimensional data and the optimal weights, which is critical to quick detection since, in such settings, there are not many post-change samples. 
Simulation results and real data examples are provided to validate the good performance of the proposed method.
\end{abstract}

\section{Introduction}

Sequential change detection is a classic problem in statistics and information theory. The goal is to detect the change in the underlying distribution as quickly as possible after it occurs. There is a wide range of applications, including sensor networks \cite{VenuSensor2010}, seismology \cite{change-point-seismology07}, social networks \cite{raginsky_OCP}, power systems \cite{chen2015quickest}, and genomics \cite{siegmund2011detecting}. Many classic results and detection procedures have been developed, see \cite{poor-hadj-QCD-book-2008,detectAbruptChange93,tartakovsky2014sequential}. However, many widely used methods assume a parametric form of the distributions before and after the change. For high-dimensional data, such parametric methods can be difficult to implement since the post-change distribution is typically unknown. We cannot have a large number of samples to estimate the distribution. Recently, there have been many interests in developing a non-parametric change detection procedure for high-dimensional streaming data when we have limited post-change samples. 

We focus on a type of {\it distribution-free} methods based on empirical distributions. Compared with parametric methods, such non-parametric tests are more flexible and can be more applicable for various real-world situations. They tend to perform better when (i) the data does not follow a parametric distribution or (ii) we do not have enough samples to estimate the underlying distribution reliably. However, one particular challenge is to establish performance guarantees and improve the sample efficiency of the non-parametric test statistic \cite{ramdas2015decreasing}.  

In this paper, we develop a new data-driven distribution-free sequential change detection procedure based on the {\it weighted} $\ell_2$ divergence between empirical distributions as the test statistic, which is related to the idea of {\it testing closeness} between two distributions \cite{batu2000testing}. We introduce ``weights'' that are design parameters, which can be particularly important in achieving good performance in practice when we do not have a large number of samples. We show the optimality of the proposed procedure in achieving the theoretical \textcolor{black}{lower bound of the sample complexity required for a low-risk test that meets the specifications} and characterize the proposed detection procedure's theoretical performance based on weighted $\ell_2$ divergence for data in a sliding window. Moreover, we develop practical optimization procedures for selecting the optimal weights and the low-dimensional projections for high-dimensional data. The good performance of the proposed procedure is demonstrated using synthetic and real data.

The rest of the paper is organized as follows. Section~\ref{sec:pre} introduces preliminaries about the problem set-up and reviews related work. Section~\ref{sec:ell2} proposes the weighted $\ell_2$ test and shows its optimality in $\ell_2$ sense. Section~\ref{sec:change} studies sequential change detection using the proposed statistic.  Section~\ref{sec:optimization} discusses different aspects to optimize the parameters involved in the proposed method. Section~\ref{sec:simulation} and  Section~\ref{sec:realdata} demonstrate the performance of the proposed detection procedure using simulation and real-data study. All proofs are delegated to appendices.

\section{Preliminaries}\label{sec:pre}

\subsection{Weighted $\ell_2$ Divergence for Two-Sample Test}\label{sec:}

We start by considering the problem of testing closeness between two discrete distributions from samples observed.  The problem set-up is as follows. Let $\Omega$ be an $n$-element {\it observation space}. For notational convenience, we identify $\Omega$ with $\{1,...,n\}$ and a probability distribution on $\Omega$ with a vector $p\in\Delta_n=\{p\in\bR^n: p\geq0,\sum_i p_i=1\}$, where $p_i$ is the probability mass of the $i$-th element in $\Omega$.
Suppose we are given two independent sample sets:
$$
X^1=\{x^{1}_{1}, \ldots, x^{1}_{n_1}\} \;\;\mbox{and}\;\; X^2=\{x^{2}_{1}, \ldots, x^{2}_{n_2}\}, 
$$
where $x^{1}_{1},\ldots,x^{1}_{n_1}\overset{\text{iid}}{\sim}p$, and $x^{2}_{1},\ldots,x^{2}_{n_2}\overset{\text{iid}}{\sim}q$. Our goal is to design a {\it test} which, given observations $X^1$ and $X^2$, claims one of the following hypotheses:
\[
H_0: \, p=q,   \quad H_1: \,  \|p-q\|_2\geq \epsilon\|p\|_2, %\|p-q\|_2\geq \epsilon,
\]
where $\|\cdot\|_2$ is the $\ell_2$ norm on $\bR^n$ and $\epsilon>0$ is a parameter that represents the relative difference of magnitude. Note that the alternative hypothesis $H_1$ considered here is slightly different but closely related to the traditional setting where $H_1$ is defined as $\|p-q\|_2\geq \epsilon$. 
%We aim to emphasize the relative difference, and this is of particular interest in quality control where the goal is to detect when the process is out of control and calibrate the threshold based on the nominal distribution $p$. 

Define the type-I risk of a test as the probability of rejecting hypothesis $H_0$ when it is true, i.e., the probability of claiming $\|p-q\|_2\geq \epsilon\|p\|_2$ when $p = q$. The type-II risk is the probability of claiming $p=q$ when $\|p-q\|_2\geq \epsilon\|p\|_2$. 
We aim at building a test which, given $0<\alpha,\beta<1/2$, has the type-I risk at most $\alpha$ (which we call {\em at level $\alpha$}), and the type-II risk at most $\beta$ ({\em of power $1-\beta$}); and we aim to meet these specifications with sample sizes $n_1$ and $n_2$ as small as possible.

We propose a new type of test by considering a family of distance-based divergence between empirical distributions of the two sets of observations. More specifically, we consider tests that reject the null hypothesis $H_0$ (and accept the alternative $H_1$) when
\[
D(X^1, X^2) > \ell,
\]
where $D(\cdot,\cdot)$ is a proxy for the {weighted} $\ell_2$ divergence between distributions $p$ and $q$ underlying $X^1$ and $X^2$, and $\ell$ is a data-dependent (random) threshold. 

Our motivation for considering the $\ell_2$ divergence for the non-parametric test is twofold. First, the $\ell_2$ divergence-based test has a certain (near) optimality that we will show in Section~\ref{sec:ell2}. Second, the $\ell_2$ divergence is more {\it robust} compared to other divergences such as $\chi^2$-divergence, which is commonly used when $\ell_1$ separation between distributions is of interest. The $\chi^2$-divergence becomes numerically difficult to evaluate when there are ``small'' $p_i$ (meaning some atoms have small probability), while the $\ell_2$ distance remains bounded and in such cases. Similar argument holds for the $\alpha$-- and $\beta$--divergences \cite{chernoff1952measure,renyi1961measures,basu1998robust,cichocki2010families,amari2016information} and detection statistic for robust change detection \cite{zhang2019robust}. Moreover, here we focus on a new {\it weighted} $\ell_2$ divergence, which emphasizes atoms that contribute most to $\|p-q\|_2$. 

Below, we derive the optimality in the sense that the sample size needed by the $\ell_2$ divergence test matches (up to a moderate constant factor) the minimax lower bound in \cite{chan2014optimal}.  More specifically, the optimal sample complexity of testing two distributions $p$ and $q$ with  $\|p-q\|_2\geq \epsilon/\sqrt{n}$ is $O(1)\sqrt{n}/ \epsilon^2$ (i.e., it requires $O(1)\sqrt{n}/ \epsilon^2$ samples from each distribution), whence any low-risk test requiring $\epsilon$ to be on the order of $O(1)/\sqrt{n}$ requires a sample size at least $O(1)\sqrt{n}$.

We also propose a general framework for selecting optimal test parameters utilizing convex optimization. Namely, when extra prior information about unknown distributions $p$ and $q$ is available, such information can be used to improve the test's quality by efficiently adjusting weights, which are adaptive to the ``closeness'' at different parts of the two distributions. Furthermore, we develop optimal projection for dimension reduction of high-dimensional data by maximizing the Wasserstein distance between two samples \cite{mueller2015principal}. This sheds light on the potential extension of the proposed test statistic to high-dimensional continuous distributions.

\subsection{Related Work}\label{sec:review}

There is a long history of studying similar problems in both statistics and computer science. In statistics, a two-sample test is a fundamental problem in which one aims to decide if two sets of observations are drawn from the same distribution \cite{lehmann2006testing}, with a wide range of applications \cite{larsen1976statistics}. Available approaches to the two-sample test can be largely divided into two categories: parametric and non-parametric. The parametric approach assumes that the data distribution belongs to certain parametric families, but the parameters can be unknown \cite{cramir1946mathematical}. The non-parametric setting does not impose any assumption on the underlying distribution and therefore is widely applicable to real scenarios. 

Classical approaches focus on the so-called ``goodness-of-fit'' test to decide whether the observations follow a pre-specified distribution. Non-parametric goodness-of-fit tests can be generalized for two-sample (and multi-sample) tests; in this case, the focus is the asymptotic analysis when the sample size goes to infinity. For instance, the Kolmogorov-Smirnov test \cite{smirnov1939estimate}, and the Anderson-Darling test \cite{Hajek1999} focus on univariate distributions and compute divergences between the empirical cumulative distributions of two (and multi) samples. The Wilcoxon-Mann-Whitney test \cite{mann1947test,kruskal1952use} is based on the data ranks and is also limited to univariate distributions. Van der Waerden tests are based on asymptotic approximation using quantiles of the standard Gaussian distribution \cite{conover1980practical,ksample}. The nearest neighbors test for multivariate data is based on the proportion of neighbors belonging to the same sample \cite{schilling1986multivariate}. 

There is much work aimed at extending univariate tests to the multivariate setting. A distribution-free generalization of the Smirnov two-sample test was proposed in \cite{bickel1969distribution} by conditioning on the empirical distribution functions. Wald-Wolfowitz run test and Smirnov two-sample test were generalized to multivariate setting using minimal spanning trees in \cite{friedman1979multivariate}. A class of distribution-free multivariate tests based on nearest neighbors was studied in
\cite{bickel1983sums, henze1988multivariate, schilling1986multivariate}, and a multivariate $k$-sample test based on Euclidean distance between sample elements was proposed in \cite{szekely2004testing}. Some recent work includes methods based on maximum mean discrepancy (MMD) \cite{gretton2012kernel} and the Wasserstein distance \cite{ramdas2017wasserstein}. 
In particular, the $\ell_2$ test enables us to draw a conclusion directly based on comparing empirical distributions. Compared with existing methods such as the MMD test, which requires a huge gram matrix when the sample size is large, the $\ell_2$ test enables us to choose weights flexibly to better serve the testing task.

Another line of research in theoretical computer science deals with {\em closeness testing}. 
It was first studied in \cite{batu2000testing,batu2013testing}, in which the testing algorithm with sub-linear sample complexity was presented; the lower bound to the sample complexity was gave in \cite{valiant2011testing}; a test that meets the optimal sample complexity was proposed in \cite{chan2014optimal}; see \cite{rubinfeld2012taming} and \cite{canonne2015survey} for recent surveys. The $\ell_2$ case has also been studied in  \cite{batu2000testing, goldreich2011testing, chan2014optimal}, and optimal algorithms are given. Many variants of closeness testing have also been studied recently. In \cite{acharya2014sublinear}, sublinear algorithms were provided for generalized closeness testing. In \cite{bhattacharya2015testing}, the closeness testing was studied under the case where sample sizes are unequal for two distributions. In  \cite{diakonikolas2017near}, a nearly optimal algorithm for closeness testing for discrete histograms was given. In \cite{acharya2018differentially}, the problem was studied from a differentially private setting.

Outstanding early contributions of sequential change detection mainly focus on parametric methods \cite{page1954,page1955test,Shiryaev1963,Lorden1971} and is well-summarized in recent books \cite{lai2001sequential,tartakovsky2014sequential}. Recently, there have been growing interests in the non-parametric hypothesis test used in change detection problems. In \cite{boracchi2018quanttree}, the ``QuantTree'' framework was proposed to define the bins in high-dimensional cases recursively, and the resulted histograms are used for change detection. In \cite{chen2019sequential}, a sequential change detection procedure using nearest neighbors was proposed. \textcolor{black}{In the seminal work \cite{lau2018binning}, a binning strategy was developed to discretize the sample space to construct the detection statistic to approximate the well-known generalized likelihood ratio test. The binned detection statistic's asymptotic properties were studied, and it was shown to be asymptotically optimal when the pre-and post-change distributions are discrete. Note that here we do not rely on likelihood ratios and assume the pre- and post-change distributions are unknown, and all we have are some possible ``training data.''}

\section{Weighted $\ell_2$ divergence test}\label{sec:ell2}

Our goal in this section is to develop a test statistic, the weighted $\ell_2$ divergence, used as the basic building block of the change detection procedure. We aim to construct a test with the following properties. When applied to two independent sets of size $N$,  i.i.d. samples $\{x_1^1,...,x_N^1\}$ and $\{x_1^2,...,x_N^2\}$ drawn from unknown distributions $p,q\in\Delta_n$, the test \begin{enumerate}
\item[(i)] rejects the null hypothesis with probability at most a given $\alpha$ under $H_0:p=q$;
\item[(ii)] accepts the null hypothesis with probability at most a given $\beta$ when there is a \textcolor{black}{relative} difference ``of magnitude at least a given $\epsilon>0$,'' i.e., under $H_1:\|p-q\|_2\geq \epsilon\|p\|_2$.
\end{enumerate}
We want to meet these reliability specifications with as small sample size $N$. %Moreover, we will show that the {\it weighted $\ell_2$ divergence test} we will build has certain optimality.

\subsection{Test Statistic}\label{sec:ell2-stat}

The main ingredient of weighted $\ell_2$ divergence test is the {\it individual test} built as follows. Let us fix ``weights'' $\sigma_i\geq0$, $i=1,...,n$, and let $\Sigma=\Diag\{\sigma_1,...,\sigma_n\}$ be \textcolor{black}{a diagonal matrix with diagonal entries being $\sigma_1,\ldots,\sigma_n$.} 
Given $\{x_1^1,...,x_N^1\}$ and $\{x_1^2,...,x_N^2\}$, we divide them into two {\it consecutive} \textcolor{black}{(left)} parts $E,E'$, of cardinality $L$ each, and \textcolor{black}{(right)} parts $F,F'$, of cardinality $R$ each, respectively. \textcolor{black}{Note that the cardinality $L$ and $R$ are at most $N/2$ and can be less than $N/2$ if we do not use all $N$ samples.} Set
\begin{equation}\label{eq:para}
\gamma={R\over L+R},\,\bgamma=1-\gamma={L\over L+R},\,M={2LR\over L+R}=2\gamma L=2\bgamma R.
\end{equation}
Let $\omega,\omega',\zeta,\zeta' \in \Delta_n$ be the empirical distributions of observations in sets $E,E',F,F'$, and $\chi$ be the weighted $\ell_2$ test statistics defined as
\begin{equation}
\chi=(\omega-\zeta)^T\Sigma(\omega'-\zeta') = \sum_{i=1}^n\sigma_i(\omega_i-\zeta_i)(\omega'_i-\zeta'_i).
\label{chi-def}
\end{equation}
The {\it weighted $\ell_2$ divergence test} $\cT$ claims a change if and only if
\begin{equation*}\label{event}
|\chi|>\ell,
\end{equation*}
where $\ell$ is the threshold. The following lemma summarizes the properties of $\cT$:
\begin{proposition}[Test Properties]\label{thm:concentration} Let $\cT$ be the weighted $\ell_2$ divergence test applied to a pair of samples drawn from distributions $p,q \in \Delta_n$, and let the threshold $\ell$ satisfy
\begin{equation}\label{ell}
\begin{array}{c}
\ell\geq 2\sqrt{2}\theta M^{-1}\sqrt{\sum_i\sigma_i^2p_i^2},
\end{array}
\end{equation}
for some  $\theta\geq1$. Then
\begin{enumerate}
\item Risk: The type-I risk of $\cT$ is at most $1/\theta^2$;
\item Power: Under the assumption
\begin{equation}\label{magnitude}
\begin{array}{rcl}
\sum_i\sigma_i(p_i-q_i)^2&>& \ell+ 2\sqrt{2}\theta \bigg[M^{-1/2}\sqrt{\sum_i\sigma_i^2(p_i-q_i)^2(\gamma p_i+\bgamma q_i)}+ \\
&&\hspace{55pt}M^{-1}\sqrt{\gamma\sum_i\sigma_i^2p_i^2+\bgamma\sum_i\sigma_i^2q_i^2}\bigg],
\end{array}
\end{equation}
the power of $\cT$ is at least $1-3/\theta^2$.
\end{enumerate}
\end{proposition}
For simplicity, in the rest of this section we assume that $\sigma_i=1$, $1\leq i\leq n$, so the left hand side of \eqref{magnitude} reduces to $\left\Vert p - q\right\Vert_2^2$. In Section \ref{sec:optweight}, we will discuss how to utilize the non-uniform weights.

\subsection{Special Case: $\ell_2$ Test With Uniform Weights}

The individual test $\cT$ in the previous section has two drawbacks: (i) to control the type-I risk, the threshold $\ell$ in \eqref{ell} specifying $\cT$  must be chosen with respect to the magnitude $\|p\|_2$ which is typically  {\it unknown}; (ii) to achieve a small type-I risk of $\cT$ we need to set a large $\theta$, thus resulting in poor power of the test. This section will show that we can reduce these limitations by ``moderately'' increasing the sample sizes. To simplify the notation, from now on, we use the fixed value $\theta=3$ (i.e., the type-I risk is at most $1/9$ and the power is at least $2/3$), and use $M=L=R$ as a special case of the definition in \eqref{eq:para}. 

The testing procedure will be as follows. We first give the Algorithm~\ref{alg} to specify the threshold $\ell$ that satisfies the condition \eqref{ell} with high probability and then introduce the testing procedure. 

\subsubsection{Specifying threshold}

When the nominal distribution $p$ is {\it unknown}, we perform a {\it training-step} -- use part of the first set of observations  to build, with desired reliability $1-\delta$, a tight upper bound $\varrho$ (\textcolor{black}{the output of Algorithm~\ref{alg}}) on the squared norm $\|p\|_2^2$ of the unknown distribution $p$ such that 
\begin{equation}\label{target}
\bP\left[ \|p\|_2^2\leq\varrho\leq 3\|p\|_2^2 \right] \geq 1-\delta,
\end{equation}
where the probability is taken with respect to the observations sampled from distribution $p$. 
\par
The training-step is organized in Algorithm~\ref{alg:pilot}, where the input parameter $S$ is defined as
\begin{equation}\label{eq:S}
S:= \min\left\{S \in \mathbb{N}: \sum_{k=S}^{2S} {2S \choose k} \left(\frac13\right)^k \left(\frac23\right)^{2S-k} \leq \frac{\delta}{\lceil \log_2(n)\rceil} \right\}.
\end{equation}
The definition in \eqref{eq:S} has an intuitive explanation: $S$ is the smallest number such that in $2S$ independent tosses of a coin, with probability of getting a head in each toss being $\leq 1/3$, the probability of getting at least $S$ heads does not exceed $\delta/m$, where $m=\lceil \log_2(n)\rceil$.
 
\begin{algorithm}
\SetAlgoLined
{\bf Input:} Samples $X^1:=\{x_1,\ldots,x_{N}\}$; Reliability $1-\delta$; $m=\lceil \log_2(n)\rceil $; $S$ in \eqref{eq:S}\;
{\bf Output:} A tight upper bound $\varrho$ on $\|p\|_2^2$ satisfying the condition \eqref{target}\;
\For{$i = 1,\ldots,m$}{
$\rho_i=2^{-i/2}$\;
Set $P_i\in\mathbb R_{+}$ as the solution to
$$
3\left[2^{7/4}P_i^{-1/2}\rho_i^{3/2}+2P_i^{-1}\rho_i\right]={1\over 3}\rho_i^2;
$$\\
Set $Q_i=\lceil P_i\rceil$\;
{Use $4S$ consecutive segments, of cardinality $Q_i$ each, of the sample $X^1$ to build $2S$ pairs $\{(\xi_s,\xi^\prime_s), s=1,\ldots,2S\}$ of empirical distributions}\;
Set $\theta_s=\xi_s^T\xi_s^\prime$ for $s=1,\ldots,2S$\;
Set $\Theta_i$ as the median of $\theta_1,...,\theta_{2S}$\;
\If{$\Theta_i\geq2\rho_i^2/3$ or $i=m$ or $N=|X^1|<4SQ_{i+1}$ (running out of sample)}{
Terminate. }
}
{\bf Return} $\varrho = \Theta_i+\rho_i^2/3$.
\caption{Training-step to estimate a tight upper bound on $\|p\|_2^2$.}\label{alg:pilot}
\label{alg}
\end{algorithm}

Properties of the training-step in Algorithm~\ref{alg:pilot} can be summarized as follows:
\begin{proposition}[Bounding $\|p\|_2^2$]\label{prop:estimate} Let $\rho_i=2^{-i/2}$ and $i(p)$ be the smallest $i\leq m$ such that $\rho_i\leq \|p\|_2$ (note that $i(p)$ is well defined due to $\rho_m\leq n^{-1/2}$). Assume that the size of the first group of sample $X^1$ is at least $4SQ_{i(p)}$. Then the probability
for the training-step to terminate in the first $i(p)$ stages and to output $\varrho$ satisfying the condition \eqref{target} is at least $1-\delta$, where $\delta$ is the reliability tolerance specifying the training-step. Besides this, the number of observations utilized in a successful training-step is at most
\begin{equation}\label{nini}
4SQ_{i(p)}=O(1)\ln(\ln(n)/\delta)/\|p\|_2.
\end{equation}
\end{proposition}

\subsubsection{Testing procedure} After $\varrho$ is built,  we use the part of the first sample $X^1$ {\it not} used in the training-step and the entire second sample $X^2$ to run $K=K(\alpha,\beta)$
individual tests to make a decision. Here $\alpha< 1/2$ and $\beta<1/2$ are pre-specified upper bounds on the type-I and type-II risks of the testing problem, and $K(\alpha,\beta)$ is the smallest integer such that the probability of getting at least $K/2$ heads in $K$ independent tosses of a coin is
\begin{itemize}
\item[(i)]
$\leq\alpha$, when the probability of getting head in a single toss is $\leq 1/9$,
 \item[(ii)]
 $\geq 1-\beta$, when the probability of getting head in a single toss is $\geq 2/3$.
\end{itemize}
It is easy to check that $K\leq O(1)[\ln(1/\alpha)+\ln(1/\beta)]$. 

The $k$-th individual test is applied to two $2M$-long segments of observations taken first from the sample $X^1$ (and these are non-overlapping with the training-step observations), and second from $X^2$, with non-overlapping segments of observations used in different individual tests. Here the positive integer $M$, same as the reliability tolerances $\delta$, $\alpha$, $\beta$, is a parameter of our construction, and the threshold $\ell$ for individual tests is chosen as
\begin{equation}\label{resultingell}
\ell=6\sqrt{2}M^{-1}\sqrt{\varrho}.
\end{equation}
After running $K$ individual tests, we claim $H_1$ if and only if the number of tests where $H_1$ is claimed is at least $K/2$.
The properties of the resulting $\ell_2$ test are presented as follows:
\begin{theorem}[Sample Complexity]\label{atmostbeta}
Consider the $\ell_2$ test above with design parameters $\delta$, $\alpha,\beta\in(0,1/2)$ and $M$. Then for properly selected absolute constants $O(1)$, the following holds true. Let $p,q$ be the true distributions from which $X^1$ and $X^2$ are sampled, and let the size $N$ of $X^1,X^2$ satisfies
\begin{equation}\label{condI}
N\geq O(1)\left[\ln(\ln(n)/\delta)/\|p\|_2+[\ln(1/\alpha)+\ln(1/\beta)]M\right].
\end{equation}
Then
\begin{enumerate}
\item The probability for the training-step in Algorithm~\ref{alg:pilot} to be successful is at least $1-\delta$, and when it happens there are enough observations to carry out $K$ subsequent individual tests.
\item Under the condition that the training-step is successful:
%(which happens with probability at least $1-\delta$),
\begin{enumerate}
\item The type-I risk (claiming $H_1$ when $p=q$) is at most $\alpha$;
\item For every $\epsilon>0$, with positive integer $M$ satisfying
\begin{equation}\label{Mislarge}
M\geq O(1){1\over\epsilon^2\|p\|_2},
\end{equation}
the type-II risk (claiming $H_0$ when $\|p-q\|_2\geq\epsilon\|p\|_2$)
% \begin{equation}\label{epsilonchange}
% \|p-q\|_2\geq\epsilon\|p\|_2,
% \end{equation}
is at most $\beta$.
\end{enumerate}
\end{enumerate}
\end{theorem}

\subsection{Near-Optimality of Proposed $\ell_2$ Divergence Test}

From the above analysis, when testing a difference of magnitude $\|p-q\|_2\geq\epsilon\|p\|_2$, reliable detection  is guaranteed when the size $N$ of samples $X^1$ and $X^2$ is at least $O(n^{1/2}\epsilon^{-2})$ (due to the fact that $\|p\|_2\geq n^{-1/2}$), with just logarithmic in the reliability parameters factors hidden in $O(\cdot)$. We will show that the $O(n^{1/2})$ sample size is the best rate can achieve unless additional {\it a priori} information on $p$ and $q$ is available. %The precise statement is as follows:
\begin{proposition}[Optimality]\label{thm:optimality}
Given cardinality $n$ of the set $\Omega$ and sample size $N$. For  i.i.d. $N$-observation samples $X^1$ and $X^2$, suppose there exists a low-risk test that can detect reliably any difference of magnitude $\|p-q\|_2\geq\epsilon\|p\|_2$ for $0<\epsilon<1/2$ such that
\begin{enumerate}
\item for every distribution $p$, the type-I risk is at most a given $\alpha<1/2$, and
\item for every distributions $p,q$ satisfying $\|p-q\|_2\geq \epsilon\|p\|_2$, the type-II risk is at most a given $\beta<1/2$.
\end{enumerate}
Then $N\geq O(1)\sqrt{n}$, with a positive absolute constant $O(1)$ that depends on $\alpha,\beta,\epsilon$.
\end{proposition}

\subsection{Illustrating Example: Quasi-Uniform Distribution} \label{sec:quasi}

Now we present an illustrative example using ``quasi-uniform'' distributions. Assume that the nominal distribution $p$ and the alternative distribution $q$ are {\it quasi-uniform}, i.e., there exists a known constant $\kappa$ satisfying $2\leq\kappa\leq n$ such that 
$\|p\|_\infty\leq\kappa/n$ and $\|q\|_\infty\leq\kappa/n$. Since $\|x\|_2\leq\sqrt{\|x\|_1\|x\|_\infty}$, we have
$
\max[\|p\|_2,\|q\|_2]\leq \sqrt{\kappa/n},
$
and hence the threshold
\begin{equation}\label{qunifa}
\ell=6\sqrt{2}M^{-1} \sqrt{\kappa/n}.
\end{equation}
satisfies the condition \eqref{ell} with $\theta=3$ (recall that we are in the case of uniform weights $\sigma_i\equiv 1$).
With this choice of $\ell$, the right hand side of condition \eqref{magnitude} is at most $6\sqrt{2} [ 2M^{-1} \sqrt{\kappa/n} + M^{-1/2}\sqrt{\kappa/n} \|p-q\|_2 ]$. 
To ensure the validity of condition (\ref{magnitude}) with $\theta=3$, it suffices to have 
\[
\begin{array}{c}
\|p-q\|_2^2\geq 6\sqrt{2} \left[ 2M^{-1} \sqrt{\kappa/n} + M^{-1/2}\sqrt{\kappa/n} \|p-q\|_2 \right],
\end{array}
\]
which holds when 
\begin{equation}\label{qinfb}
\begin{array}{c}
\|p-q\|_2^2\geq O(1)M^{-1} \sqrt{\kappa/n}, 
\end{array}
\end{equation}
with properly selected moderate absolute constant $O(1)$. For quasi-uniform distributions, $\|p-q\|_2$ is no larger than $2\sqrt{\kappa/n}$. Therefore, for $
\|p-q\|_2\geq \lambda n^{-1/2}
$ with some $\lambda\in(0,2\sqrt{\kappa}]$, the sample size $M$ should satisfy 
\begin{equation*}\label{M}
M\geq O(1){\sqrt{\kappa n}\over\lambda^2}
\end{equation*} 
in order to ensure condition (\ref{qinfb}). We see that in the case of $L=R$, given $\alpha\ll1$, $\beta\ll1$, the sample size of
$$
O(1)[\ln(1/\alpha)+\ln(1/\beta)]{\sqrt{\kappa n}\over\lambda^2}
$$
ensures that for the  $\ell_2$ test with the threshold (\ref{qunifa}), its type-I risk and type-II risk are upper bounded by $\alpha$ and $\beta$, respectively.

\begin{figure}[h!]
\centerline{
\includegraphics[width=0.4\textwidth]{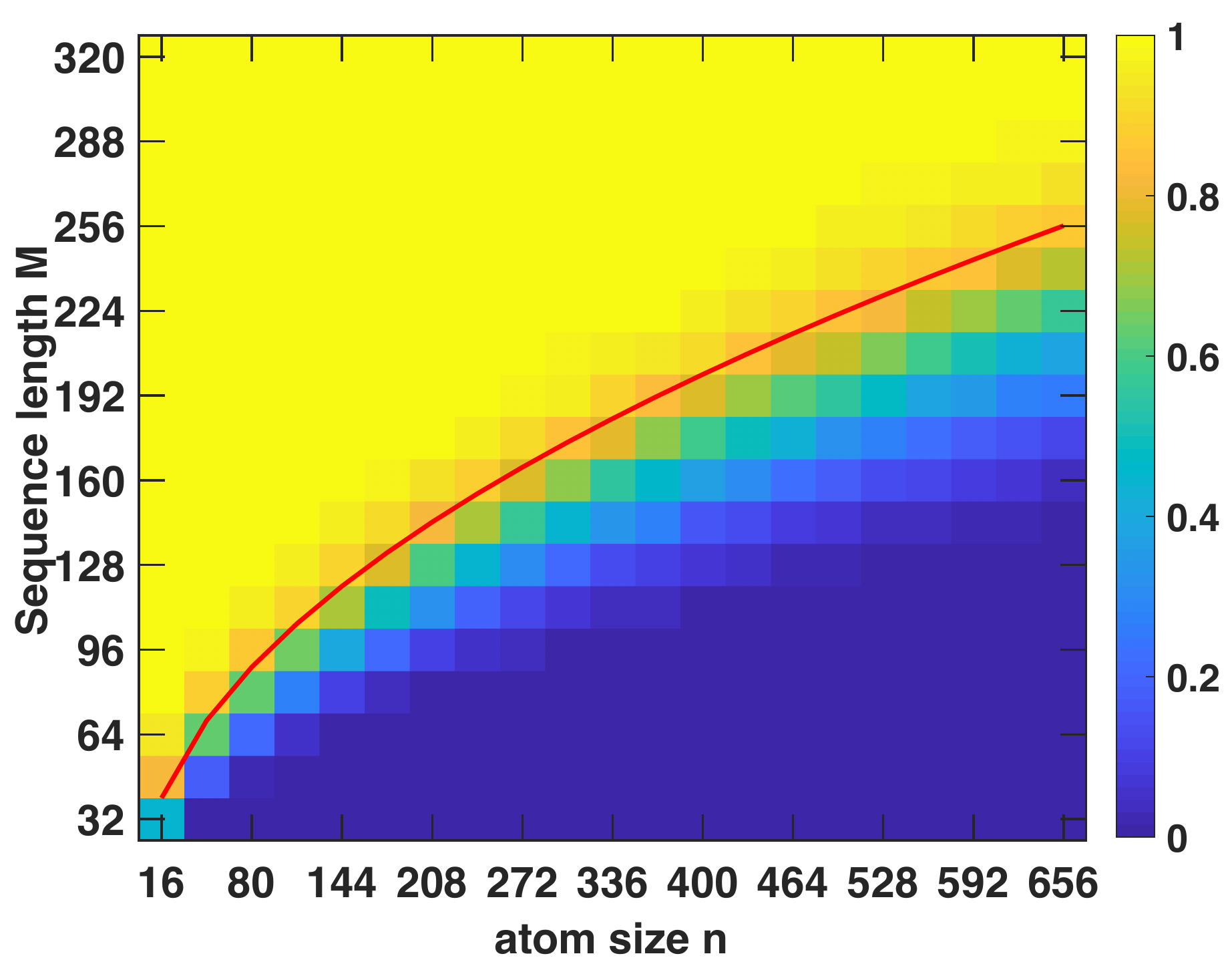} }
\caption{Validation of the theoretical $O(\sqrt{n})$ bound by plotting the empirical test power  of ``quasi-uniform'' in Section \ref{sec:quasi}, averaged over 1000 random trials. The Type-I risk is controlled to be less than $10^{-3}$. The theoretical lower bound to sample complexity $O(\sqrt{n})$ is shown in red line, which match the empirical phase-transition ``watershed''.}
\label{fig:simulation1}
\vspace{-0.3in}
\end{figure}

In the following, we provide numerical examples to validate the optimality results in Proposition~\ref{thm:optimality}. Suppose the support size $n$ is even and set $L=R=M$ for simplicity. The experiment set-up is described as the following two steps:
\begin{enumerate}[label=(\roman*)]
\item Draw two $n/2$-element subsets independently, $\Omega_1$ and $\Omega_2$, of $\Omega$ from the uniform distribution on the family of all subsets of $\Omega$ of cardinality $n/2$. 
\item The samples $X^1$ are i.i.d. drawn from the uniform distribution on $\Omega_1$, denoted as $p$; and the second group of samples $X^2$ are i.i.d. drawn from the uniform distribution on $\Omega_2$, denoted by $q$. 
\end{enumerate}
Therefore we have $\max[ \left\| p \right\|_2 ,\left\| q \right\|_2] \leq \sqrt{2/n}$, implying that we can set the threshold as
\[\ell = 12M^{-1}n^{-1/2}.\]
In all simulations, the individual test was applied. We perform the simulation for various $n$ and $M$ values. The power is shown in Fig. \ref{fig:simulation1}, averaged over 1000 trials. The results show that for magnitude $\left\|p-q\right\|_2 = O(1/\sqrt{n})$, at least $O(\sqrt{n})$ samples are required in order to detect the difference between $p$ and $q$ with high probability.

\subsection{Comparison With Standard Identity Test}

The most attractive feature of the $\ell_2$ test is that to eliably test a distribution separation of relative ``magnitude'' $\epsilon$, the sizes of samples should be of on the order of $O(\sqrt{n}/\epsilon^2)$, just logarithmic in $n$ and reliability tolerance (type-I and type-II risks) factors hidden in $O(\cdot)$. In contrast, the conceptually simplest approach to test -- comparing norm $\|\widehat{p}-\widehat{q}\|_\pi$  of the difference of two empirical distributions with theoretical quantile of the recovery error as measured in the same norm -- turns out to require sample sizes of the order of $n$. To the best of our knowledge, the only test proposed in the literature where the required sample size grows with $n$ sublinearly is the ``identity test'' from \cite{valiant2017automatic} designed for the test {\it with known in advance nominal distribution $p$ (the null distribution).} It is shown in \cite{valiant2017automatic} that the number of observations allowing for the identity test to detect reliably the shift $p$ to $q$ is, up to logarithmic in the reliability parameters factors, $\cM_p={\sqrt{n}/\|p-q\|_1^2}$. When applying in the same situation of known $p$ the $\ell_2$ test, invoking our Theorem \ref{atmostbeta}, the required number of observations is $\cM=\max[1,{\|p\|_2/\|p-q\|_2^2}]$.

A natural way to compare the quality of the identity test and the $\ell_2$ tests is to look at the {\it performance ratio}
$
\cC={\cM_p/\cM};
$
$\cC\gg 1$ means that the $\ell_2$  test outperforms the identity test, while $\cC\ll1$ indicates the opposite situation.
In the case of ${\|p\|_2/\|p-q\|_2^2}\geq1$, the performance ratio becomes
$$
\cC={\sqrt{n}\over\|p\|_2}\cdot{\|p-q\|_2^2\over\|p-q\|_1^2},
$$
which is the product of two factors in the range of $[O(\sqrt{n}),O(n)]$ and $[O(1/n),1]$, respectively. Thus, $\cC$ is in the range $[O(1/\sqrt{n}),O(n)]$, implying that no one of the two tests always outperforms the other one. However, {\it with proper implementation}, $\ell_2$ test always outperforms the identity test.

Specifically, with $p$ known in advance, let us apply $\ell_2$ test to an {\it equivalent} problem built as follows. 
Given $p$, let us define $k_i=\lfloor p_i n\rfloor$ for each point  $i\in\Omega=\{1,...,n\}$ and define its {\it $p$-children} $(i,k)$, $1\leq k\leq k_i$. Let $\overline{\Omega}$ be the set of children of all points from $\Omega$. Observing a random variable $\xi$ taking values in $\Omega$, we can convert these observations into observations of random variable $\overline{\xi}$ taking values in $\overline{\Omega}$ as follows: observing value $i$ of the random variable $\xi$, we throw a perfect dice with $k_i$ faces and take $(i,k)$ as the value of $\overline{\xi}$, where $k$ being the observed face of the dice. Note that the distribution $\overline{p}$ of $\overline{\xi}$ is obtained from the distribution $p$ of $\xi$ by spreading, for every $i\in\Omega$, the probability mass $p_i$ of $i$ in equal portions $p_i/k_i$ between the children of $i$. In particular, $\overline{p}$ has all entries at most $1/n$, while the cardinality $\overline{n}$ of $\overline{\Omega}$ clearly is at most $2n$.
Now, given samples $X^1$ and $X^2$ of cardinality $N$ each, the distribution being (known) $p$ and (unknown) $q$, let us use the above randomization to convert the samples into i.i.d. samples, of the same size, drawn from distributions $\overline{p}$ and $\overline{q}$, and use these resulting samples to perform testing. With this approach, the sample size allowing to detect reliably the difference between $p$ and $q$, up to factor logarithmic in reliability tolerances,  will become $\overline{\cM}=\max[1,\|\overline{p}\|_2/ \|\overline{q}-\overline{p}\|_2^2]$.
In this case of $\|\overline{p}\|_2/\|\overline{q}-\overline{p}\|_2^2\geq 1$, we have 
\begin{equation}\label{eeee}
{\cM_p\over\overline{\cM}}={\sqrt{n}\over\|\overline{p}\|_2}\cdot{\|\overline{p}-\overline{q}\|_2^2\over \|p-q\|_1^2}.
\end{equation}
Taking into account that $\|\overline{p}\|_2\leq \sqrt{2/n}$, $\|p-q\|_1=\|\overline{p}-\overline{q}\|_1$ and $\overline{n}\leq 2n$, the first factor in the right hand side of (\ref{eeee}) is at least $O(n)$, and the second is at least $O(1/n)$, meaning that $\cM_p/\overline{\cM}$ is at least $O(1)$. Therefore, with the above implementation, $\ell_2$ test is never much worse and in fact can be much better than the identity test.

\subsection{Comparison With Generalized $\pi$-tests} \label{sec:ell_pi}

Here we compare $\ell_2$ test with {\it $\pi$-tests} based on estimating $\|\cdot\|_\pi$-distance between distribution $p$ and $q$ given samples, and will make an argument that $\ell_2$ metric is preferred than other metrics.  Specifically, given samples $X^1$ and $X^2$ of cardinality $2N$ each and $\pi\in[1,\infty]$, consider the test  as follows:
\begin{enumerate}[label=(\roman*)]
\item We split the first set of sample $X^1$ into two non-overlapping samples -- {\it training subsample}  $x^\tr$ and {\it testing subsample} $x^\ts$, of cardinality $N$ each, and compute the empirical distributions $\hat{f}[x^\tr]$, $\hat{f}[x^\ts]$ of these samples. In the same fashion, we divide the second set of sample $X^2$, thus getting two $N$-observation non-overlapping training subsample $\bar{x}^\tr$ and testing subsample $\bar{x}^\ts$ and the corresponding empirical distributions $\hat{f}[\bar{x}^\tr]$, $\hat{f}[\bar{x}^\ts]$;
\item We use the full training sample $[x^\tr,\bar{x}^\tr]$ to build a {\it separator} vector $e\in\bR^n$ such that
$$
e^T[\hat{f}[x^\tr]-\hat{f}[\bar{x}^\tr]]=\|\hat{f}[x^\tr]-\hat{f}[\bar{x}^\tr]\|_\pi, \quad  \|e\|_{\pi_*}=1,\quad \pi_*={\pi\over\pi-1},
$$
and set $h=e/\|e\|_\infty$;
\item We compute the quantity $\chi=h^T[\hat{f}[x^\ts]-\hat{f}[\bar{x}^\ts]]$ and claim $H_1$ if and only if
$$
\chi>\ell,
$$
where the threshold $\ell$ is the parameter of our test.
\end{enumerate}
The properties of the resulting test can be summarized as follows. Let $\rho_{\beta,n,\pi}(N)$ be an upper bound on the $(1-\beta)$-quantile of the $\|\cdot\|_\pi$-error of recovering a probability distribution $p$ on observation space $\Omega=\{1,...,n\}$ by the empirical distribution of $N$-element i.i.d. sample drawn from $p$. In other words, for every probability distribution $f$ on $\Omega$, the empirical distribution $\hat{f}[x^N]$ associated with drawn from $f$ $N$-element i.i.d. sample $x^N$ satisfies the relation
$$
\bP\{\|\hat{f}[x^N]-f\|_\pi>\rho_{\beta,n,\pi}(N)\}\leq\beta.
$$
It can be shown that a tight, within absolute constant factor, choice of $\rho_{\beta,n,\pi}(N)$ is given by
\begin{equation}\label{quantile}
\rho_{\beta,n,\pi}(N)=2\sqrt{3\ln(2n/\beta)}n^\varkappa/\sqrt{N},\,\varkappa=\max\left[{1\over\pi}-{1\over 2},0\right].
\end{equation}
We have the following result:
\begin{theorem}[Properties of Generalized $\pi$-test]\label{ell1_risk}
Let $N$ be a positive integer and $\pi\in[1,\infty]$. For  $\alpha\in(0,1)$,  let $\tau_{N}(\alpha)$ be the positive root of the equation
\begin{equation}\label{equation}
{\tau^2\over 1+{2\over 3}\tau}={2\over N}\ln(2/\alpha),
\end{equation}
and define $\tau_{N}'(\alpha) = \min\{1,\tau_{N}(\alpha)\}$.
For the $\pi$-test utilizing samples $X^1$ and $X^2$ of cardinality $2N$ each, if the threshold is set as
\begin{equation}\label{ellthresh}
\ell = 2\tau_{N}'(\alpha)
\end{equation}
(let us denote this test $\cT_{\pi,\alpha}^N$), one has:
\begin{enumerate}
\item The type-I risk of $\cT_{\pi,\alpha}^N$ is $\leq \alpha$;
\item Given $\beta\in(0,1/4)$, under the assumption
\begin{equation}\label{final}
\|p-q\|_\pi\geq 6\rho_{\beta,n,\pi}(N)+2\tau_{N}'(\alpha),
\end{equation}
where $p$, $q$ are true distribuiton underlying samples $X^1$ and $X^2$,
the type-II risk of the test $\cT_{\pi,\alpha}^N$ is at most $4\beta$.
\end{enumerate}
\end{theorem}
Note that the power of $\cT_{\pi,\alpha}^N$ is the same as of the ``straightforward'' test, which uses samples $X^1$ and $X^2$ to build empirical approximations of the corresponding distributions and claims $H_1$ when the $\|\cdot\|_\pi$-distance between these approximations is ``essentially larger'' than (theoretical upper bound on the)  $(1-\alpha)$-quantile of the $\|\cdot\|_\pi$-norm of the approximation error. The advantage of the $\pi$-test over the straightforward one is that the $\pi$-test is less conservative in terms of the type-I risk. 
%
%Note that the power of $\cT_{\pi,\alpha}^N$ is the same as of the ``straightforward'' test, which uses samples $X^1$ and $X^2$ to build empirical approximations of the corresponding distributions and claims $H_1$ when the $\|\cdot\|_\pi$-distance between these approximations is ``essentially larger'' than (theoretical upper bound on the)  $(1-\alpha)$-quantile of the $\|\cdot\|_\pi$-norm of the approximation error. The advantage of our $\pi$-test as compared to the straightforward one is that the $\pi$-test is less conservative as far as the type-I risk is concerned. Indeed, what matters in this respect for $\pi$-test is the theoretical upper bound on $(1-\alpha)$-quantile of the value of a specific linear form of the approximation error, and not a bound  on the $(1-\alpha)$-quantile of $\|\cdot\|_\pi$-norm of this error (this is what matters for the straightforward test), and standard bounds in the first of these two bounding problems seemingly are essentially less conservative than in the second one.
%
Let us compare $\pi$-tests with the $\ell_2$-test. It is easy to see that there is no reason to use $\pi$-tests with $\pi>2$; indeed, when $\pi>2$, the right hand side in (\ref{final}) is the same as when $\pi=2$, see (\ref{quantile}), while the left hand side in (\ref{final}) is non-increasing in $\pi$. Thus, considering the power of $\pi$-tests, we can restrict $\pi\in[1,2]$.

Next, assume that $\alpha$ and $N$ are such that $\tau=\tau_N(\alpha)<1$ (this is the only nontrivial case, since otherwise the test $\cT_{\pi,\alpha}^N$ never claims $H_1$). In this case, (\ref{equation}) implies that $\sqrt{\Theta}\leq\tau\leq {4\over3}\sqrt{\Theta}$, with $\Theta=2\ln(2/\alpha)/N$. Finally, assume that $n/\beta\geq 2/\alpha$; under these assumptions (\ref{quantile}) says that the first term in the right hand side of (\ref{final}) dominates the second one, so that (\ref{final}) implies that in order for $H_1$ to be detected with reliability $1-4\beta$, we should have
\begin{equation}\label{basically}
\|p-q\|_\pi\geq O(1)n^{{1\over\pi}-{1\over 2}}\sqrt{\ln(2n/\beta)/N},
\end{equation}
with moderate absolute constant $O(1)$. Condition (\ref{basically}) implies that the best of $\pi$-tests is the one with $\pi=2$. Indeed, when $\pi\in[1,2]$, by H\"older inequality, we have $\|p-q\|_\pi\leq n^{{1\over \pi}-{1\over 2}}\|p-q\|_2$, implying that whenever condition (\ref{basically}) is satisfied by some $p,q$ and some $\pi\in[1,2]$, it is satisfied by the same $p,q$ with $\pi=2$.

Another implication is that the $\ell_2$ test significantly outperforms the best of $\pi$-tests. Indeed, by Theorem \ref{atmostbeta}, for all $\delta,\alpha,\beta\in(0,1/2)$, samples of cardinality $N$,
with
$$
N\geq N_{\ell_2}=O(1)\left[\epsilon^2\ln(\ln(n)/\delta)+\ln(1/\alpha)+\ln(1/\beta)\right]{1\over\epsilon^2\|p\|_2}
$$
allow to ensure conditional, the probability of the condition being at least $1-\delta$, type-I and type-II risks are at most $\alpha$ and $\beta$, respectively.
By (\ref{basically}), to ensure similar reliability properties for $\pi$-test with $\pi=2$, we should have a much larger, for small $\|p\|_2$, the sample size
$$
N\geq N_2=O(1){\ln(2n/\beta)\over\epsilon^2\|p\|_2^2}.
$$
In the ``extreme case'' of $\|p\|_2=O(1)/\sqrt{n}$,
$
N_2/N_{\ell_2}
$ is as large as $O(1)\sqrt{n}$. The discussion above supports the near-optimality of the $\ell_2$ test.

\section{Sequential change detection procedures}\label{sec:change}

In this section, we construct the change detection procedure based on the proposed weighted $\ell_2$ divergence test. Change detection is an important instance of the sequential hypothesis test, but it has unique characteristics that require a separate study due to different performance metrics considered. Since we do not know the change location, we have to perform scanning when forming the detection statistic. We discuss two settings: the offline scenario where we have fixed samples and the online setting where the data come sequentially. 

\subsection{Offline Change Detection by ``Scan'' Statistic} 

In the offline setting, we observe samples $X^T=\{x_1,...,x_{T}\}$ on a time horizon $t=1,...,T$, with $x_t$'s taking values in an $n$-element set $\Omega=\{1,\ldots,n\}$. Assume there exists time $K\in\{1,...,T\}$ such that for $t\leq K$, $x_t$ are i.i.d. drawn from some
{\it pre-change} distribution $p$, and for $t\geq K+1$, $x_t$ are i.i.d. drawn from the {\it post-change} distribution $q$. Our goal is to design a test which, based on the samples $X^T$, decides on the null hypothesis $K=T$ (``no change")  versus the alternative $K<T$ (``change''). Meanwhile, we want to control the probability of false alarm to be at most a given $\alpha > 0$,
and under this restriction to make the probability of successfully detecting the change as large as possible, at least when $K$ and $T-K$ both are moderately large and $q$ ``significantly differs"
from $p$.

We use the proposed test in Section~\ref{sec:ell2} to construct a scan statistic for change detection. Given $T$, we select a collection of {\it bases} $B_j$, $1\leq j\leq J$. A base $B$ is a segment of $\{1,...,T\}$ partitioned into three {\it consecutive} parts: {\it pre-change} part $B_\lf$, middle part $B_\md$, and {\it post-change} part $B_\rg$; the last instant in $B_\lf$ is the first instant in $B_\md$, and the first instant in $B_\rg$ is by 1 larger than the last instant in $B_\md$. For example: 
$
B_\lf=\{1,...,10\},\,B_\md=\{10,11\},B_\rg=\{12,...,20\}$. 
We associate with base $B$ an {\it individual test} $\cT_B$ which operates with observations $\{x_t, t\in B_\lf\cup B_\rg\}$ only. This test aims at deciding on two hypotheses: (1) ``No change:'' there is no change on $B$, that is, either $K$ is less than the first, or larger than or equal to the last time instant from $B$; (2) ``Change:'' the change point $K$ belongs to the middle set $B_\md$.

Given $\alpha>0$ and a base $B$, we call individual test $\cT_B$ associated with this base {\it $\alpha$-feasible},  if the probability of false alarm for $\cT_B$ is at most $\alpha$, meaning that whenever there is no change on the base $B$ of the test, the probability for the test to claim change is at most $\alpha$. 
Our ``overall'' test $\cT$ works as follows: we equip bases $B_j$, $1\leq j\le J$, with tolerances $\alpha_j>0$ such that
% \begin{equation}\label{sumeps}
% \begin{array}{c}
$\sum_j\alpha_j=\alpha$, 
% \end{array}
% \end{equation}
and then associate with each base $B_j$ with a $\alpha_j$-feasible individual test $\cT_{B_j}$ (as given by the $\ell_2$ test in Section~\ref{sec:ell2-stat}). Given observations $X^T$, we perform one by one the individual tests in some fixed order, until either (i) the current individual test claims change; when it happens, the overall test claims change and terminates, or (ii) all $J$ individual tests are performed and no one of them claimed change; in this case, the overall test claims no change and terminates.

\begin{proposition}[False Alarm Rate for Offline Change Detection]
With the outlined structure of the overall test and under condition $\sum_j\alpha_j=\alpha$,
the probability of false alarms for $\cT$ (of claiming change when $K=T$) is at most $\alpha$.
\end{proposition}

%Our primary goal in the sequel is to understand how to associate with given base $B$ and tolerance $\alpha>0$ a $\alpha$-feasible test $\cT_B$ with ``presumably high power,'' meaning that whenever the change point $K$ belongs to $B_\md$ and the (properly defined) magnitude of the change $q-p$ is not small, the probability of miss -- of $\cT$ claiming ``no change'' -- is small. There is a tradeoff among a few factors: the probability of miss, the magnitude of change, and the cardinalities of pre- and post-change parts of $B$. % and we intend to arrive at nearly as good tradeoff as possible. 
%After we known how to design feasible tests given their bases, we intend to address the (conceptually simpler) question of how to design these bases.

\subsection{Online Change Detection} 

Instead of giving a fixed duration of samples in the offline setting, the observations arrive sequentially for online detection tasks. The goal is to detect the change as quickly as possible, under the constraint that the false alarm rate is under control. 

\begin{figure}[H]
\vspace{-0.1in}
\centerline{
\includegraphics[width = 0.6\textwidth]{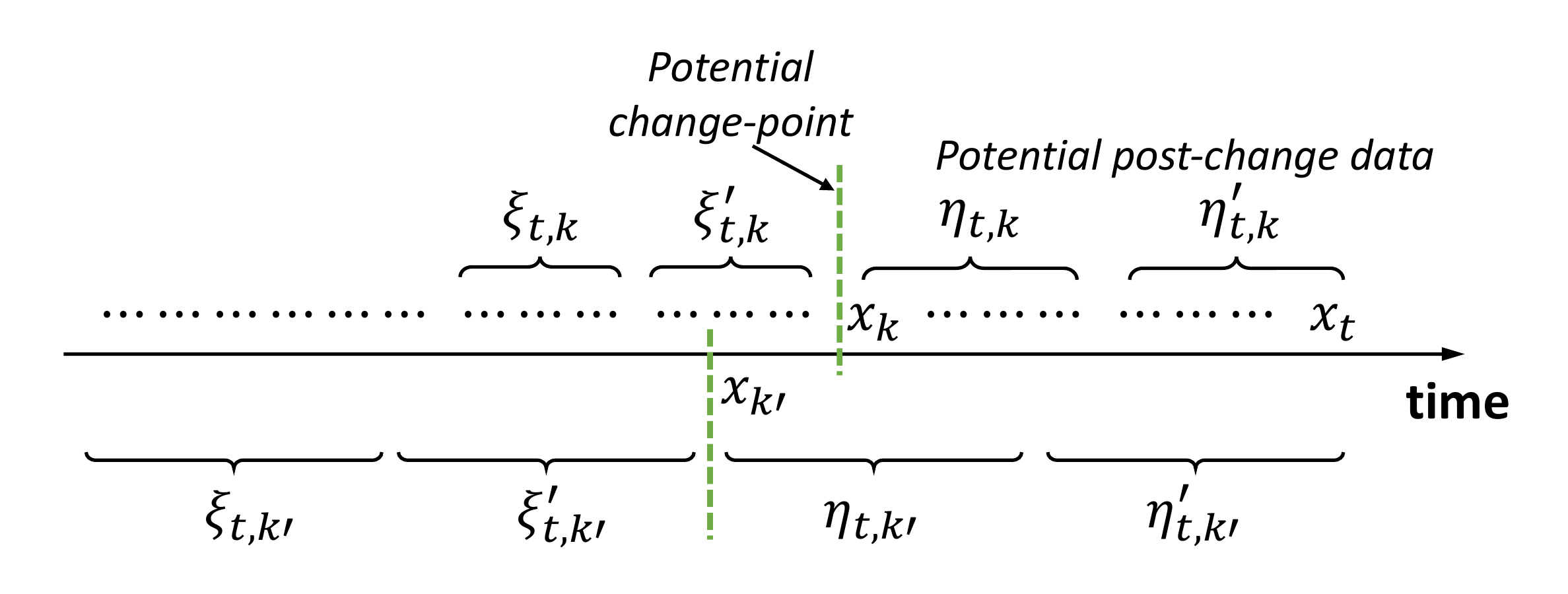} }
\caption{Illustration of the sequential change detection procedure.}
\label{fig:slidingwindow}
\vspace{-0.1in}
\end{figure}

The proposed detection procedure based on $\ell_2$ test is illustrated in Fig.~\ref{fig:slidingwindow}. Given a sequence $\{x_t,t=1,2,\ldots\}$, as each time $t$, we search over all possible change-points $k<t$. In particular, we form two sequences before $k$ and two sequences between $[k,t]$ with the length all equal to $M_{t,k} = \lceil (t-k)/2 \rceil$; their corresponding empirical distributions are denoted as $\xi_{t,k}$, $\xi_{t,k}'$, and $\eta_{t,k}$, $\eta_{t,k}'$. The detection statistic $\chi_{t,k}$ is formed as:
\begin{equation}\label{eq:detect_stat_online}
\begin{array}{c}
\chi_{t,k} = M_{t,k}(\xi_{t,k}-\eta_{t,k})^\intercal \Sigma (\xi'_{t,k}-\eta'_{t,k}).
\end{array}  
\end{equation}
We note that the multiplicative term $M_{t,k}$ can be viewed as a {\it scaling} parameter (which is proportional to the standard deviation of the test statistic) such that the variance of $\chi_{t,k}$ is of a constant order as $t-k$ increases. 

The online change-point detection procedure is given by a stopping time 
\begin{equation}\label{eq:stoptime}
\cT:= \inf\{t: \max_{0\leq k\leq t} \chi_{t,k} \geq b\},
\end{equation}
where $b$ is a pre-specified threshold that needs to be determined by controlling the false alarm rate. An intuitive interpretation of $\cT$ is that at each time $t$, we search over all possible change-pints $k<t$, and raise alarm if the maximum statistic exceeds the threshold. 
\begin{remark}[Window-Limited Procedure]
In practice, we can also adopt a {\it window-limited} version of $\cT$ as
\begin{equation}\label{eq:stoptime_window}
\cT' = \inf\{t: \max_{m_0\leq t-k\leq m_1} \chi_{t,k} \geq b\},
\end{equation}
where $m_0$ and $m_1$ are the lower and upper bounds of the window size that we would like to scan for the possible changes. \textcolor{black}{Usually $m_0$ can be set such that the resulted sequences are long enough to have meaningful empirical distributions for constructing the detection statistic in \eqref{eq:detect_stat_online}. For practical considerations, we usually require the window size $m_1$ to be at least the expected detection delay, as discussed in \cite{lai1998information} where the original theoretical study of the window-limited test was proposed.} 
%satisfy $\log m_1 = o(\log(\mbox{\emph{ARL}}))$ and meanwhile let $m_1$ to be roughly the same order as the desired detection delay, as discussed in 
\end{remark}
\begin{remark}[Comparison with the Binning Approach \cite{lau2018binning}]
\textcolor{black}{We note that the binning approach in \cite{lau2018binning} also considers discretized space and scans over all possible change-points when approximating log-likelihood ratio statistics. Compared with \cite{lau2018binning} which assumes the pre-change distribution is known, the detection procedure \eqref{eq:stoptime} and its window-limited version \eqref{eq:stoptime_window} do not need the prior knowledge of the pre-change distribution. We did not use the log-likelihood ratio statistic here but scan over all neighboring time windows directly to detect any significant difference in empirical distributions.}
\end{remark}

\subsection{Theoretical Analysis for Online Change Detection}

Now we characterize the two fundamental performance metrics for sequential change detection, namely the average run length (ARL) and the expected detection delay (EDD). We cannot use the previous method in Proposition~\ref{thm:concentration} to determine the threshold because the bound is too conservative and will be intractable when the ARL is large. Here we present an asymptotic theoretical approximation. 

\subsubsection{Theoretical approximation to ARL}\label{sec:ARL}

To compute the ARL, we need to quantify the distribution of $\cT$ when data are sampled from the same distribution $p$. Intuitively, the detection statistic $\chi_{t,k}$ is small when the samples are from the same distribution. \textcolor{black}{A relatively standard result is that when the threshold $b$ tends to infinity, the stopping time's asymptotic distribution is approximately exponential when there is no change. This is proven true in various scenarios \cite{siegmund1995using,siegmund2008detecting,yakir2009multi}. The main idea is to show that the number of boundary cross events for detection statistics over disjoint intervals converges to Poisson random variable in the total variation norm; the result can be established by invoking the Poisson limit theorem for dependent samples. Detailed proofs by adapting those techniques into the specific $\ell_2$ test setting are left for future work.} 
%A commonly used approach is to approximate $\cT'$ (and $\cT$) using an exponential distribution \cite{siegmund1995using}. 
Under such approximation, we have
\[
\bP_\infty(\cT'\geq m) = \bP_\infty\left( \max_{1\leq t \leq m}\max_{m_0\leq t-k\leq m_1} \chi_{t,k} \geq b\right)  \approx e^{-\lambda m},
\]
where $\bP_\infty$ is the probability measure when the change-point equals to $\infty$, i.e., the change never happens; and $\bE_\infty$ denotes the corresponding expectation under this probability measure. Therefore, we only need to compute the probability $\bP_\infty(\cT'\geq m)$ and find the parameter $\lambda$, then the expectation of $\cT'$ equals $1/\lambda$. We adopt the change-of-measure transformation \cite{yakir2013extremes,xie2013sequential,SiegmundYakirZhang2010} and characterize the local properties of a random field. We first quantify the correlation between $\chi_{t,k}$ ane $\chi_{\tau,s}$ in order to find the probability $\bP_\infty(\cT' \geq m)$ theoretically.
\begin{proposition}[Temporal Correlation of Sequential Detection Statistics]\label{prop:corr}
Suppose all samples are i.i.d. drawn from the same distribution $p$, denote $M = t-k = \tau-s$, then the correlation between $\chi_{t,k}$ and $\chi_{\tau,s}$ is
\[
\mathrm{Corr}(\chi_{t,k},\chi_{\tau,s}) = 1 - \frac{2}{M}|t-\tau|  - \frac{2}{M}|k-s|  + \frac{(t-\tau)^2+(s-k)^2}{M^2}.
\]
\end{proposition}
Based on the correlation result, we have the following Theorem characterizing the ARL of the proposed $\ell_2$ sequential detection procedure. The main idea is to use a linear approximation for the correlation between detection statistics $\chi_{t,k}$ and $\chi_{\tau,s}$. Then the behavior of the detection procedure can be related to a random field. By leveraging the localization theorem \cite{SiegmundYakirZhang2010}, we can obtain an asymptotic approximation for ARL when the threshold $b$ is large enough (in the asymptotic sense).  Define a special function $v(\cdot)$ which is closely related to the Laplace transform of the overshoot over the boundary of a random walk \cite{siegmund2007statistics}:
\begin{equation}\label{eq:overshot}
v(x) \approx \frac{\frac{2}{x}[\Phi\big(\frac{x}{2}\big)-0.5]}{\frac{x}{2}\Phi\big(\frac{x}{2}\big)+\phi\big(\frac{x}{2}\big)}, 
\end{equation}
where $\phi(x)$ and $\Phi(x)$ are the probability density function and cumulative density function of the standard Gaussian distribution. For simplicity, we denote the variance of $\chi_{t,k}$ as 
\begin{equation}\label{eq:sigma}
\begin{array}{c}
\sigma_p^2 := \mathrm{Var}[\chi_{t,k}] = 4\left[ \sum_{i=1}^n \sigma_i^2p_i^2(1-p_i)^2 + \sum_{i\neq j}\sigma_i\sigma_j p_i^2p_j^2\right].
\end{array}
\end{equation}
\begin{theorem}[ARL Approximation]\label{thm:ARL}
For large values of threshold $b \rightarrow \infty$, the ARL of the test $\cT'$ can be approximated as
\begin{equation}
\bE_\infty[\cT'] = \frac12 b^{-1} e^{b^2/(2\sigma_p^2) }[2\pi \sigma_p^2]^{1/2} \Big / \int_{[4b^2/(m_1\sigma_p^2)]^{1/2}}^{[4b^2/(m_0\sigma_p^2)]^{1/2}} y \nu^2(y)dy(1+o(1)).  
\label{eq:ARL}
\end{equation}
\end{theorem}

The main contribution of Theorem~\ref{thm:ARL} is to provide a theoretical method to set the threshold that can avoid the Monte Carlo simulation, which could be time-consuming, especially when ARL is large. \textcolor{black}{Although there is no close-form analytical solution for $b$, when we let the right-hand side of Equation \eqref{eq:ARL} equals to a specific ARL value (lower bound), we can numerically compute the right-hand side of \eqref{eq:ARL} for any given threshold value $b$. Then we search over a grid to find the corresponding threshold values. Table \ref{tab:threshold1} validates the approximation's good accuracy by comparing the threshold obtained from Equation \eqref{eq:ARL} and compares it with that obtained by the Monte Carlo simulation.} In detail, we generate 2000 independent trials of data from nominal distribution $p$ and perform the detection procedure $\cT'$ for each trial; the ARL for each threshold $b$ is estimated by the average stopping time over 2000 trials. In Table~\ref{tab:threshold1}, we report the threshold obtained through Monte Carlo simulation (as a proxy for the ground-truth) and on the approximation \eqref{eq:ARL}, for a range of ARL values. The ARL values in Table~\ref{tab:threshold1} correspond to the lower bound of an ARL; since ARL will increase when increasing the threshold. So if we have a good approximation, this can help us to calibrate the threshold and control the false alarm rate. The results in Table~\ref{tab:threshold1} indicate that the approximation is reasonably accurate since the relative error is around 10\% for all specified ARL values. It is worth mentioning that ARL is very sensitive to the choice of threshold, making it challenging to estimate the threshold with high precision.  However, the EDD is not that sensitive to the choice of the threshold, which means that a small difference in the threshold will not significantly change EDD.
\begin{table}[h!]\setlength\tabcolsep{8pt}
\renewcommand\arraystretch{1.0}
\begin{center}
\caption{Comparison of the threshold $b$ obtained from simulations and the approximation \eqref{eq:ARL}. Scanning window $m_0=10,m_1=50$, support size $n=20$, nominal distribution $p$ is uniform.}
\label{tab:threshold1}
\begin{tabular}{ccccccc}
\specialrule{.08em}{0em}{0em}
ARL & 5k  &  10k  &  20k  & 30k &  40k & 50k\\
\hline
Simulation   &   2.0000  &  2.1127  &  2.2141  &  2.2857   & 2.3333 &   2.3750 \\
\hline
Theoretical & 1.8002  &  1.8762 &   1.9487  &  1.9897   & 2.0183  &  2.0398 \\
\specialrule{.08em}{0em}{0em}
\end{tabular}
\vspace{-0.3in}
\end{center}

\end{table}

%Obviously, the Monte Carlo simulation is universal for all change detection problems as long as the nominal distribution is known in advance or can be estimated reliably from historical data. However, the simulation-based method is time-consuming compared with the theoretical-based method. 

\subsubsection{Theoretical characterization of EDD} 

After the change occurs, we are interested in the expected detection delay, i.e., the expected number of additional samples to detect the change. There are a variety of definitions for the detection delay \cite{Lorden1971,poll-astat-1985,pelkowitz1987asymptotically,tartakovsky2014sequential}. To simplify the study of EDD, it is customary to consider a specific definition $\bE_0[\cT']$, which is the expected stopping time when the change happens at time 0 and only depends on the underlying distributions $p,q$. It is not always true that $\bE_0[\cT']$ is equivalent to the standard worst-case EDD in literature \cite{Lorden1971,poll-astat-1985}. However, since $\bE_0[\cT']$ is certainly of interest and is reasonably easy to approximate, we consider it as a surrogate here. We adopt the convention that there are certain pre-change samples $\{x_{-1},x_{-2},\ldots\}$ available before time 0, which can be regarded as reference samples.  

Note that for any $t>0$ and $k=0$, the sequences $\xi_{t,0}$ and $\xi'_{t,0}$ come from the pre-change distribution $p$ since they belong to the reference sequence $\{x_{-1},x_{-2},\ldots\}$, and the sequences $\eta_{t,0}$ and $\eta'_{t,0}$ are from the post-change distribution $q$. Therefore, the expectation of the detection statistic $\chi_{t,k}$ is $\bE[\chi_{t,k}] = \lceil (t-k)/2\rceil (p-q)^\intercal \Sigma (p-q)$, which determines the asymptotic growth rate of the detection statistic after the change. Using Wald's identity \cite{Siegmund1985}, we are able to obtain a first-order approximation for the detection delay, provided that the maximum window size $m_1$ is large enough compared to the EDD. 

\begin{theorem}[EDD Approximation]\label{thm:EDD} \textcolor{black}{Suppose  $b\rightarrow \infty$, with other parameters held fixed. If the window size $m_1$ is sufficiently large and greater than $2b/[(p-q)^\intercal \Sigma(p-q)]$}, then the expected detection delay 
\begin{equation}\label{eq:edd}
\bE_0[\cT] = \frac{b(1+o(1))}{(p-q)^\intercal \Sigma (p-q)/2}. 
\end{equation}
\end{theorem}

From the approximation result in Theorem~\ref{thm:ARL}, we note that the ARL is of order $O(e^{b^2/(2\sigma_p^2)})$ with respect to the threshold value $b$, which means that the threshold $b$ chosen for a fixed ARL value $\gamma$ should be on the order of $O(\log \sqrt{\gamma})$. Moreover, by taylor expansion and Csisz\'ar-Kullback-Pinsker inequality \cite{jungel2016entropy}, we have the Kullback–Leibler (KL) divergence is equivalent to squared $\ell_2$ norm up to certain constants, i.e., $C_0 \Vert p-q \Vert_2^2 \leq \mbox{KL}(p,q) \leq C_1\Vert p-q \Vert_2^2 $, therefore the EDD is also of order $b/\mbox{KL}(p,q)$, which matches the theoretical lower bound of EDD in first-order.  

\begin{remark}[Optimize weights to minimize EDD] From the EDD approximation in \eqref{eq:edd}, it is obvious that we can minimize EDD by optimizing over the weights matrix $\Sigma$. In particular, the EDD can be minimized when we can find the weights $\Sigma$ such that the weighted $\ell_2$ divergence between $p$ and $q$ is maximized. This is consistent with the subsequent discussion in Section~\ref{sec:optweight}. In particular, when we have certain prior information about the distributions $p$ and $q$, we could apply the optimization-based method in Section~\ref{sec:optweight} to find the optimal weights to reduce the detection delay.
\end{remark}

\section{Optimized weights and projection of high-dimensional data}\label{sec:optimization}

This section discusses setting optimal weights that adapt to the closeness at different elements in $\Omega$, given some {\it a prior} information on $p$ and $q$. In addition, we tackle the data high-dimensionality by adopting the Wasserstein-based principal differences analysis \cite{mueller2015principal} to find the optimal projection. %, which sheds some light on the potential application of methods for continuous distributions. 

\subsection{Optimize Weights for $\ell_2$ Test}\label{sec:optweight}

So far, we primarily focused on the case with uniform weights $\sigma_i\equiv 1$. In this section, we will discuss how to further improve performance by choosing the optimal weights. In the simplest case, when we know in advance (or can infer from additional ``training'' samples) that the distribution shift $p\to q$ (nearly) does not affect probabilities with indexes from some {\it known} set $I$, we can set $\sigma_i=0$ for $i\in I$ and $\sigma_i=1$ for $i\not\in I$. This will keep the magnitude $\sum_i\sigma_i(p_i-q_i)^2$ on the left hand side of (\ref{magnitude}), as compared to uniform weights, intact, but will  reduce the right hand side of (\ref{magnitude}).

A framework to optimize over $\sigma_i$'s is as follows. Assume that we know distributions $p,q$ belong to a set ${\cal X}\subset\Delta_n$, which is defined by a set of quadratic constraints:
\begin{equation}\label{constraints}
{\cal X} = \{p\in\Delta_n: p^TQ_k p\leq 1,\,k=1,\ldots,K\},
\end{equation}
where $Q_k\in\bR^{n\times n}$ are positive semi-definite ($Q_k\succeq0$).
\par
A natural way to measure ``magnitude of difference'' is to use $\|p-q\|_2$ (the case using $\|p-q\|_1$ can be similarly defined and solved). Assume we want to select $\sigma=[\sigma_1,...,\sigma_n]\geq0$ to make reliable detection of difference $\|p-q\|_2 \geq\rho$, for some given $\rho>0$. To achieve this, we can impose a fixed upper bound on the right hand side in \eqref{magnitude} when $p=q\in\cX$, i.e., to require $\sigma$ to satisfy
\begin{equation*}\label{semrel}
g_*(\sigma):=\max_{p\in \cX} \sum_i\sigma_i^2p_i^2 \leq a
\end{equation*}
with some given $a$, and to maximize under this constraint the quantity
$$
\begin{array}{c}
f_*(\sigma):=\min\limits_{p,q}\left\{\sum_i\sigma_i(p_i-q_i)^2:p,q\in\cX, \|p-q\|_2\geq\rho\right\}.
\end{array}
$$
For any $\sigma$ that satisfies $g_*(\sigma)\leq a$, the associated test which claims $H_1$ when the statistics (defined in \eqref{chi-def}) $|\chi| >2\sqrt{2}\theta M^{-1}\sqrt{a}$ is with type-I risk at most $1/\theta^2$. At the same time, large $f_*(\sigma)$ is in favor of good detection of distribution shift of magnitude $\|p-q\|_2\geq\rho$.
By the homogeneity in $\sigma$, we can set $a=1$ without loss of generality. 

In general, both $g_*$ and $f_*$ are difficult to compute. Therefore, we replace the problem
$$
\max_{\sigma\geq0} \left\{f_*(\sigma):g_*(\sigma)\leq1\right\}
$$
with its {\it safe tractable approximation}:
\begin{equation}\label{eq:optimal_weights}
\max_{\sigma\geq0} \left\{f(\sigma):g(\sigma)\leq1\right\},
\end{equation}
where $f$ is a concave  efficiently computable {\it lower} bound on $f_*$, and $g$ is a convex efficiently computable {\it upper} bound on $g_*$.
\par
To build $g$, note that when $p\in\cX$, the matrix $P=pp^T\in\bR^{n\times n}$ is positive semi-definite ($P\succeq0$), non-negative in each entry ($P\geq0$),  $\sum_{i,j=1}^n P_{ij}=1$, and $\Tr(PQ_k)\leq1$, $k\leq K$, by \eqref{constraints}. Consequently, the function
\begin{equation*}\label{eq:g_relax}
\begin{array}{c}
g(\sigma):=\max\left\{\Tr(\Sigma^2 P):P\succeq0,P\geq0,\sum_{i,j=1}^n P_{ij}=1,\Tr(PQ_k)\leq1,\, 1\leq k\leq K\right\}
\end{array}
\end{equation*}
with $\Sigma:=\Diag\{\sigma_1,...,\sigma_n\}$ is an efficiently computable convex upper bound on $g_*$.
Similarly, to build $f$, observe that the matrix $S=(p-q)(p-q)^T$ stemming from $p,q\in\cX$ with $\|p-q\|_2\geq\rho$ belongs to the convex set
$$
\begin{array}{c}
{\cal S}=\left\{S:S\succeq0,\sum_{i,j=1}^n|S_{ij}|\leq 4,\sum_{i,j=1}^nS_{ij}=0,\Tr(S)\geq\rho^2,\Tr(SQ_k)\leq4,\,1\leq k\leq K\right\}.
\end{array}
$$
Therefore,
$$
f_*(\sigma) \geq f(\sigma):=\min_{S\in{\cal S}}\Tr(\Sigma S)
$$
and the function $f(\sigma)$ is concave and efficiently computable.

To implement the problem \eqref{eq:optimal_weights} efficiently, we derive the tractable dual formulation in the following. Note that these constraints can be greatly simplified if $Q_k$ are {\it diagonal} matrices, especially for the high dimensional case.
\begin{proposition}[Dual Reformulation]\label{prop:dual}
The dual formulation of the optimization problem \eqref{eq:optimal_weights} is
\[
\begin{array}{rcl}
\max \quad && \lambda \rho^2 - 4\sum_k x_k  - 4\xi \\
\textrm{s.t.} \quad&& \lambda \geq0, P \succcurlyeq0, \xi \geq0, x_k \geq0, U \geq 0, W \geq 0, \Lambda \succcurlyeq 0, V \geq 0, \mu_k \geq 0, 1\leq k \leq K,\\
&& \sum_k x_kQ_k +U-W-P-rJ-\lambda I_n \succcurlyeq -\Sigma, \\
&& U_{ij} + W_{ij} \leq \xi, \quad  1\leq i\leq n, 1\leq j\leq n,\\
&& \sum_k \mu_k - \nu \leq 1, 
%& \Lambda \succcurlyeq 0, V \geq 0, \mu_k \geq 0, 1\leq k \leq K, \\
\quad - \Lambda -V + \sum_k \mu_k Q_k - \nu J \succcurlyeq \Sigma^2.
\end{array}
\]
where $\Sigma=\Diag\{\sigma_1,...,\sigma_n\}$ and $J \in \mathbb{R}^{n \times n}$ is a matrix with all elements equal to $1$.
\end{proposition}

\subsection{Comparing Optimal Weighted Versus Unweighted $\ell_2$ Divergence Test}

\begin{figure}[H]
\centerline{
\begin{tabular}{cc}
\includegraphics[width = 0.35\textwidth]{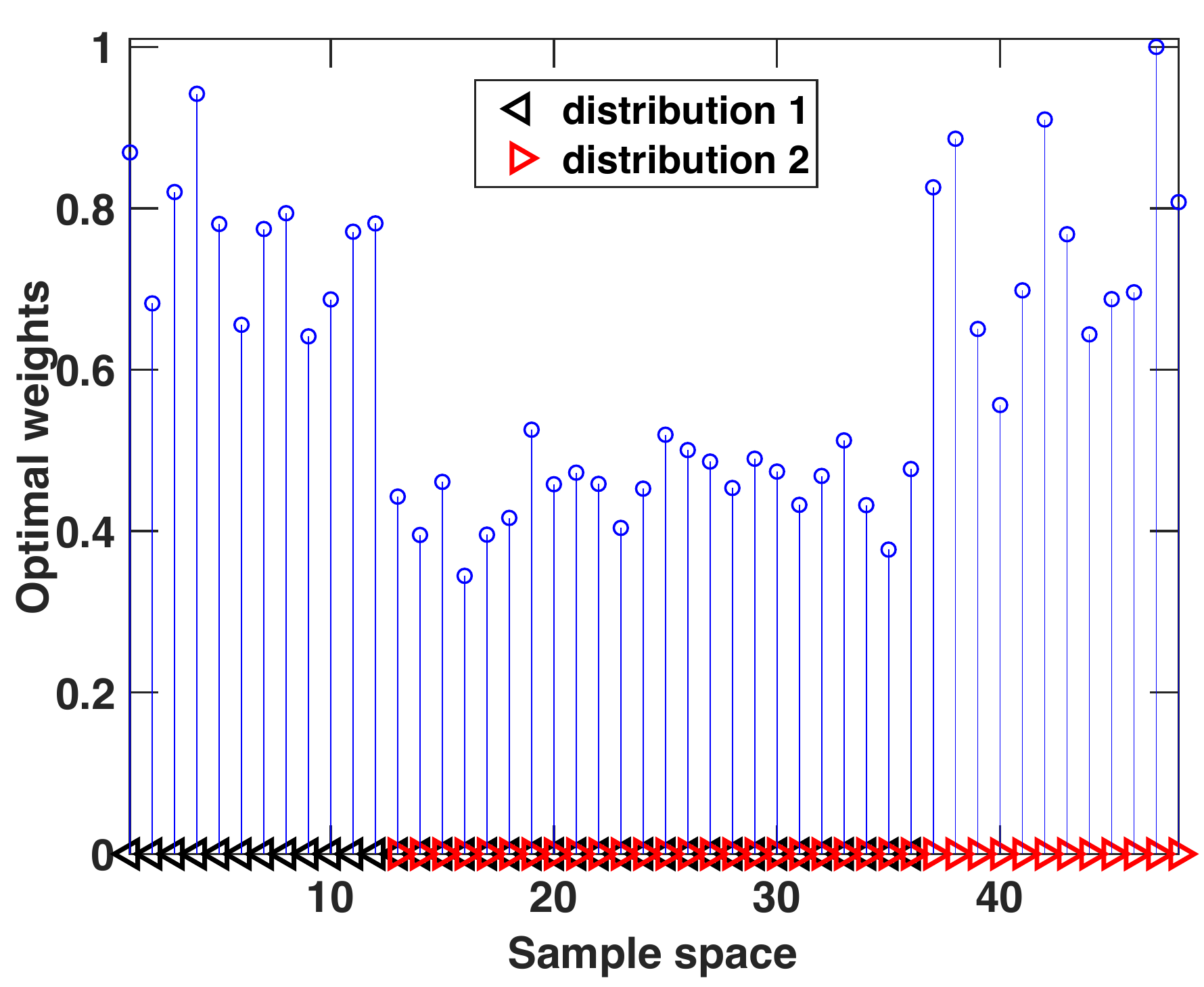} & \includegraphics[width = 0.35\textwidth]{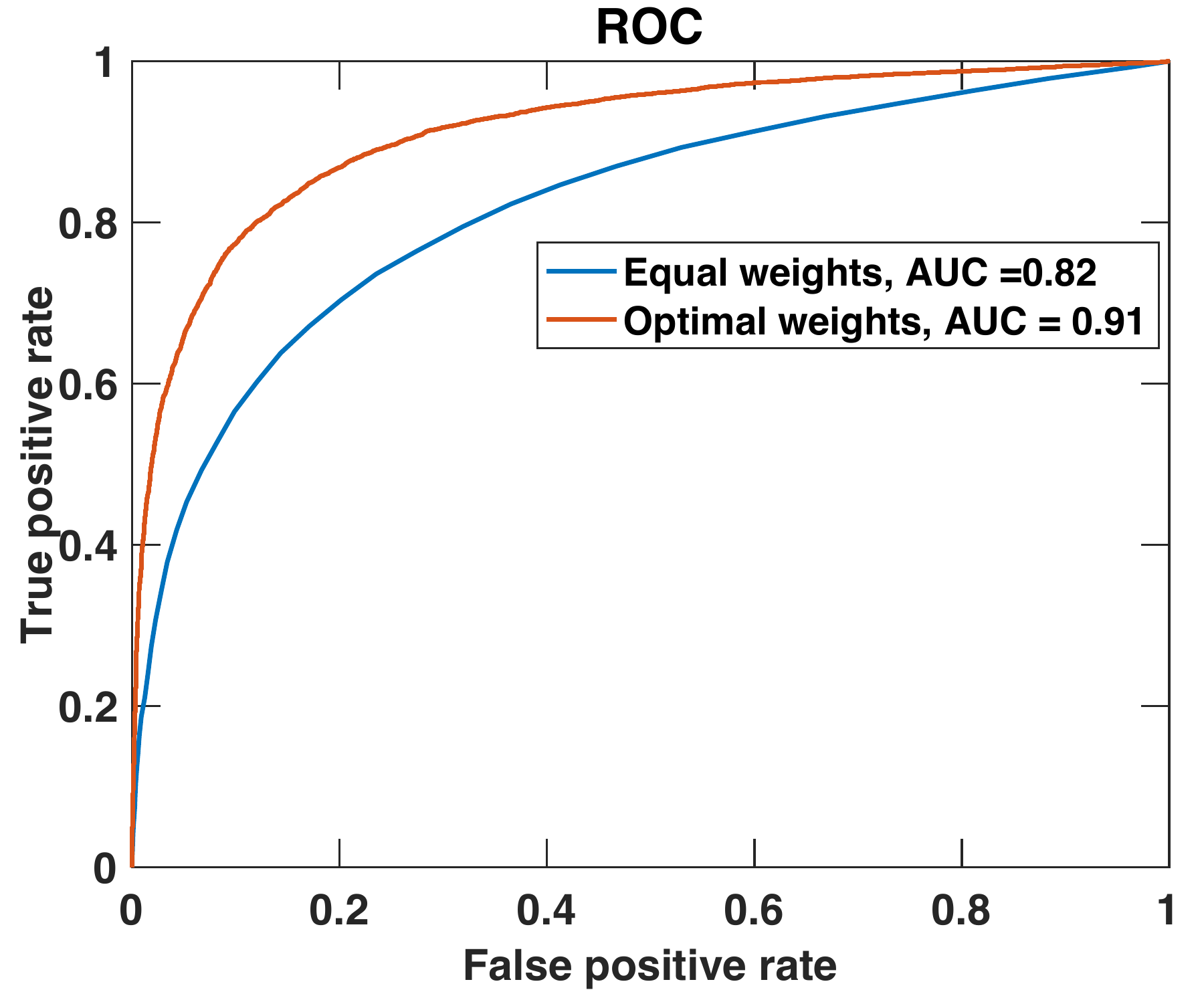} \\
(a) & (b)
\end{tabular} }
\vspace{-0.1in}
\caption{Illustration of optimal weights on a simulated example. (a): Optimal weights; (b): The ROC curves under optimal weights and equal weights.}
\vspace{-0.1in}
\label{fig:opt_weights}
\end{figure}

We present an illustrative simulation example to show the benefits of optimizing weights $\sigma$. The experimental set-up is as follows. 
Consider the sample space $\Omega=\{1,\ldots,n\}$ with $n=48$. The distributions $p$ and $q$ are set as uniform distributions on the subset $\Omega_1 = \{1,\ldots,3n/4\}$ and $\Omega_2 = \{n/4+1,\ldots,n\}$, respectively. The common support of $p$ and $q$ consists of $n/2$ elements. We first use training data to estimate the matrix $Q$ in our formulation. Specifically, we sample 32 observations from each distribution and compute the empirical distribution of all observations. This process is repeated for $m=50$ trials, and the resulting $Q$ is solved from the following optimization problem
\[
\begin{aligned}
\textrm{min}\quad & \log\det A^{-1} \\
\textrm{s.t.} \quad& \left\Vert Ap_i \right\Vert_2 \leq 1,\ i= 1,\ldots,m,
\end{aligned}
\]
where $A=Q^{1/2}$ and $p_i$ is the empirical distribution in the $i$-th trial. The volume of the ellipsoid defined with $Q$ is proportional to $\det A^{-1}$. Thus the solution to the above optimization problem is the {\it minimum ellipsoid} that contains the $m$ empirical distributions \cite{boyd2004convex}. The optimal weights are shown in Fig.~\ref{fig:opt_weights}(a). Moreover, we compare the ROC curve of the test under equal weights $\sigma_i = 1$ and optimal weights in Fig.~\ref{fig:opt_weights}(b), averaged over 10,000 trials. The result shows the benefits of using optimal weights.

\subsection{Optimal Projection for High-Dimensional Data}
\label{sec:optproj}

Assume that the two distributions $p$ and $q$, rather than discrete distributions on a given finite set, are continuous distributions on $\bR^d$. In this situation, we may try to convert observations $x\in\bR^d$ into observations $f(x)$ taking values in a {\it finite} set and apply the proposed test to the transformed observations.

One way to build $f$ is to project observations $x$ onto one-dimensional subspace and then split the range of the projection into bins. We propose to select this subspace using, when available, ``training sample'' $x_1,...,x_{2T}$, with the first $T$ observations drawn, independently
of each other, from the nominal distribution $p$, and the last $T$ observations drawn independently of each other and of $x_1,...,x_T$, from the distribution $q$. A natural selection of the one-dimensional subspace can be as follows. Denote by $e$ the unit vector spanning the subspace. Let us look at the sample empirical distributions of the projections of the observations
$x_1,...,x_{2T}$ on $e$, and try to find unit vector $e$ for which the Wasserstein distance between the distributions of the first half and the second half of the projections is as large as possible \cite{mueller2015principal}.
The distance above is, up to factor $1/T$, the quantity
\begingroup
\allowdisplaybreaks
$$
\begin{array}{rcl}
\phi(e)&=&\min\limits_{\omega_{ij, \, 1\leq i,j\leq 2T}}\left\{\sum\limits_{i,j}|e^T(x_i-x_j)|\omega_{ij}:
\begin{array}{l}
 \sum\limits_j\omega_{ij}=\begin{cases} 1,&i\leq T\\
0,&i>T\\ \end{cases} ; \sum\limits_i\omega_{ij}=\begin{cases} 0,&j\leq T\\
1,&j>T\\ \end{cases}\\
\omega_{ij}\geq0,1\leq i,j\leq 2T \\
\end{array}\right\}\\
&=&\max\limits_{\lambda}\left\{\sum_{i=1}^T\lambda_i-\sum_{i=T+1}^{2T}\lambda_i: \lambda_i-\lambda_j\leq |e^T(x_i-x_j)|,\,1\leq i,j\leq 2T\right\}\\
&=&\Phi(E[e]),
\end{array}
$$
\endgroup
where 
$$
\begin{array}{rcl}
E[e]&=&ee^T,\\
\Phi(E)&=&\max\limits_{\lambda}\left\{\sum_{i=1}^T\lambda_i-\sum_{i=T+1}^{2T}\lambda_i: \lambda_i-\lambda_j\leq \sqrt{[x_i-x_j]^TE[x_i-x_j]},\,1\leq i,j\leq 2T\right\}.
\end{array}
$$
Note that function $\Phi(E)$ is concave and the goal is to maximize $\Phi(E)$ over positive semi-definite rank-one matrices $E=ee^T$ with trace $1$. An efficiently solvable convex relaxation after relaxing the rank-one constraints is:
$$
\max_{E,\lambda}\bigg\{\sum_{i=1}^T\lambda_i-\sum_{i=T+1}^{2T}\lambda_i:
\lambda_i-\lambda_j\leq \sqrt{[x_i-x_j]^TE[x_i-x_j]}, 1\leq i,j\leq 2T,E\succeq0,\Tr(E)=1\bigg\}.
$$
After the optimal solution $E_*$ to the problem is found, we can use standard methods to obtain a reasonably good $e$, e.g., take $e$ as the leading eigenvector of $E_*$.
%or sample $N$ (e.g., 100) samples from Gaussian distribution $\mathcal N(0, E_*)$, normalize each of them to have unit $\|\cdot\|_2$-norm, thus getting $N$ unit vectors $e^t$, $t\leq N$, and select from them the one with the largest $\phi(\cdot)$. We compare it with selecting $e$ at random from the uniform distribution.

%\subsection{Comparing Optimal Projection Versus Random Projection}

\begin{figure}[H]
\centerline{
\begin{tabular}{cc}
\includegraphics[width = 0.35\textwidth]{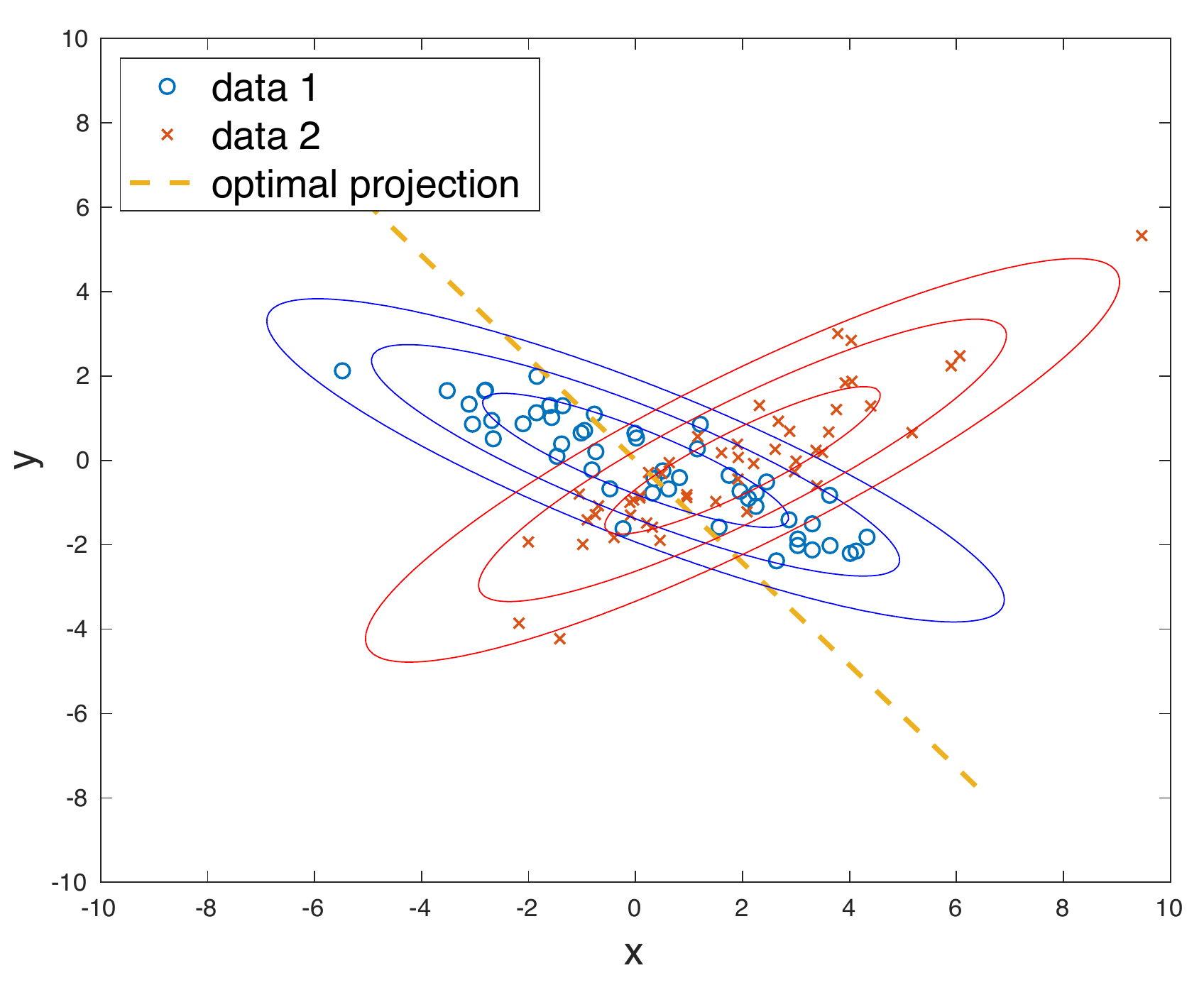}  & \includegraphics[width = 0.35\textwidth]{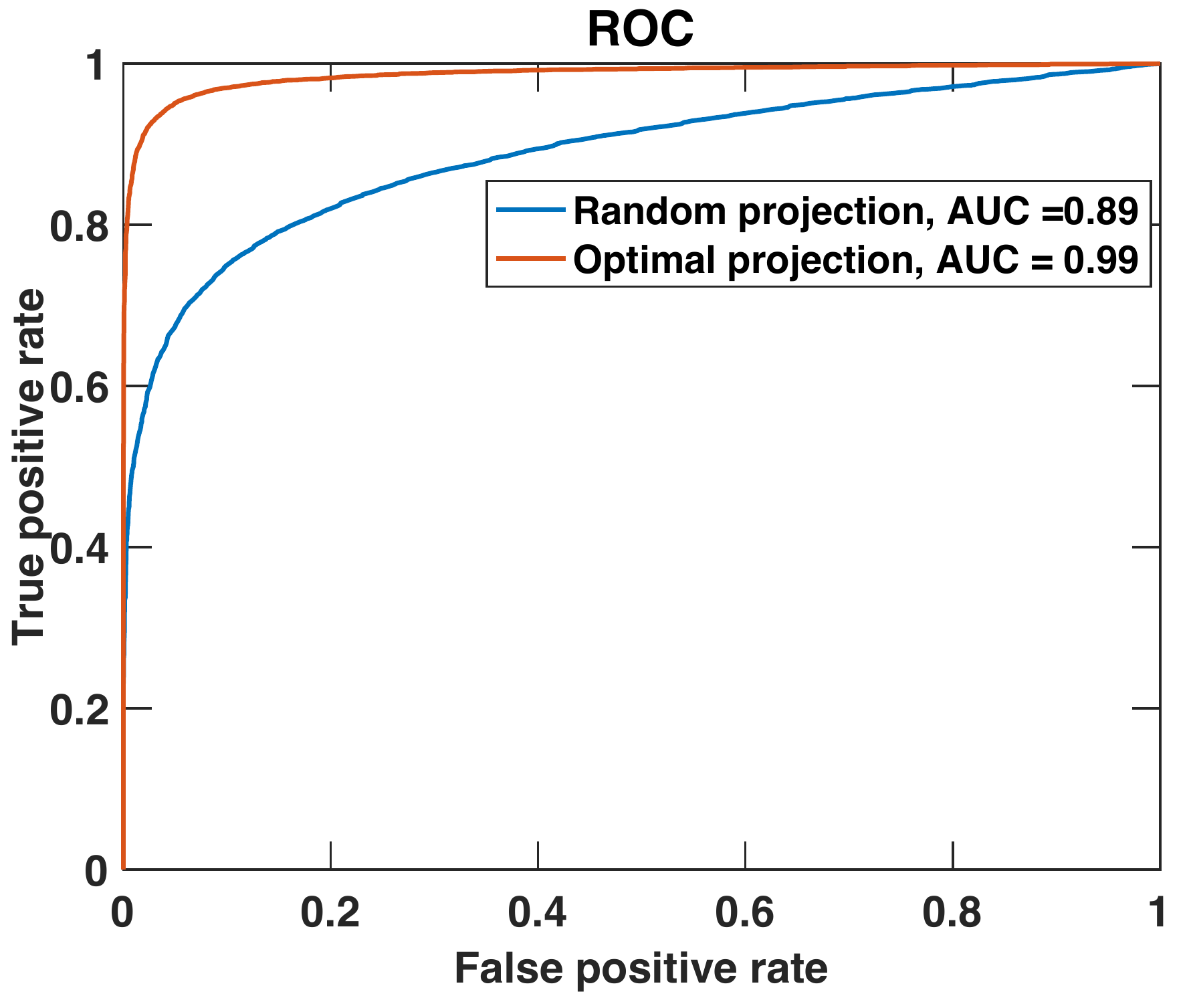}\\
(a) & (b)
\end{tabular}}
\vspace{-0.1in}
\caption{Illustration of optimal projection on simulated data. (a):  Optimal projection for two training sets; (b): The ROC curves for optimal projection and random projection.}
\label{fig:opt_proj_gaussian}
\vspace{-0.1in}
\end{figure}

Here we present a simple numerical illustration for the optimal project. Consider the two-dimensional Gaussian distributions with same mean value and different covariance structures. More specifically, let the data $X_1$ to be sampled i.i.d. from $\mathcal N(\mu_1, \Sigma_1)$ and data $X_2$ to be sampled from $\mathcal N(\mu_2, \Sigma_2)$, where
\begin{equation}\label{eq:mean_cov_shift}
\begin{array}{c}
\mu_1 = \begin{bmatrix} 0 & 0 \end{bmatrix}^\intercal, \
\mu_2 = \begin{bmatrix} 2 & 0 \end{bmatrix}^\intercal, \
\Sigma_1 = \begin{bmatrix} 5.03 & -2.41 \\ -2.41 & 1.55\end{bmatrix}, \
\Sigma_2 = \begin{bmatrix} 5.50 & 3.30 \\ 3.30 & 2.53\end{bmatrix}.
\end{array}
\end{equation}
Fig.~\ref{fig:opt_proj_gaussian} shows the optimal projection obtained from 50 training samples from each distribution (which can be seen to optimally ``separate'' the two distributions), and the ROC curve averaged over 10,000 trials that demonstrates the performance gain of the optimal projection.

\section{Numerical Examples}\label{sec:simulation}

In this section, we perform some simulations to validate the performance of the $\ell_2$ test and compare with two benchmarks: (i) the classical parametric Hotelling's $T^2$ test \cite{hotelling1947multivariate}; and (ii) the non-parametric maximum mean discrepancy (MMD) test \cite{gretton2012kernel}. More specifically, we study the test power of the two-sample test for Gaussian distributions under various dimensions. Moreover, we show the performance in change detection by studying the detection power in the offline case and the expected detection delay in the online case, respectively. 

We first introduce briefly the two benchmark procedures. 
\paragraph{Hotelling's $T^2$ statistic} The Hotelling's $T^2$ statistic is a classical parametric test designed utilizing the mean and covariance structures of data, and thus it can detect both the mean and covariance shifts \cite{hotelling1947multivariate}. Given two set of samples $\{x_1,\ldots,x_{n_1}\}$ and $\{y_1,\ldots,y_{n_2}\}$, the Hotelling's $T^2$ statistic is defined as 
\begin{equation}\label{eq:T2}
t^2 = \frac{n_1 n_2}{(n_1+ n_2)}(\bar x - \bar y)^T \widehat\Sigma^{-1} (\bar x - \bar y),
\end{equation}
where $\bar x$ and $\bar y$ are the sample mean and $\widehat\Sigma$ is the pooled covariance matrix estimate. 

\paragraph{MMD statistic} The MMD test is a non-parametric benchmark for two-sample test and change detection \cite{gretton2012kernel,li2015m}. Given a class of functions $\mathcal F$ and two distributions $p$ and $q$, the MMD distance between $p$ and $q$ is defined as 
$\mathrm{MMD}_{\mathcal F}(p,q) = \sup_{f\in\mathcal F}(\bE_{x\sim p}[f(x)] - \bE_{y\sim q}[f(y)])$. 
For MMD in reproducing kernel Hilbert spaces (RKHS), given samples $\{x_1,\ldots,x_{n_1}\}$ and $\{y_1,\ldots,y_{n_2}\}$, an unbiased estimate of squared MMD distance is given by 
\begin{equation}\label{eq:MMD}
\mathrm{MMD}_u^2 = \frac{1}{n_1(n_1-1)} \sum_{i=1}^{n_1}\sum_{j\neq i} k(x_i,x_j) + \frac{1}{n_2(n_2-1)} \sum_{i=1}^{n_2}\sum_{j\neq i} k(y_i,y_j) - \frac{2}{n_1n_2} \sum_{i=1}^{n_1}\sum_{j=1}^{n_2} k(x_i,y_j),
\end{equation}
where $k(\cdot,\cdot)$ is the kernel function associated with RKHS.

\subsection{Two-Sample Test}

Following a similar setup as in \cite{gretton2012kernel}, we investigate the performance of various tests as a function of the dimension $d$ of the sample space $\bR^d$, when both $p$ and $q$ are Gaussian distributions. We consider values of $d$ up to $256$. The type-I risk for all tests is set as $\alpha=0.05$. The sample size is chosen as $n_1 = n_2 =100$, and results are averaged over 500 independent trials. In the first case, the distributions $p,q$ have different means and the same variance. More specifically, $p = \mathcal N(0,I_d)$ and $q = \mathcal N(\mu \bm{1}/\sqrt{d}, I_d)$ with $\mu=0.8$. Note that the division of each element of the mean vector by $\sqrt{d}$ makes the difficulty of the hypothesis testing similar across all $d$ values. In the second case, the distributions $p,q$ have the same means but different variance. More specifically, $p = \mathcal N(0,I_d)$ and $q = \mathcal N(0, \Sigma)$ with $\Sigma =\Diag\{0.25,1,\ldots,1\}$, i.e., we only scale the first diagonal entry in the covariance matrix to make the hypothesis testing problem challenging to perform.

\begin{figure}[h!]
\centering
\includegraphics[width = 0.6\textwidth]{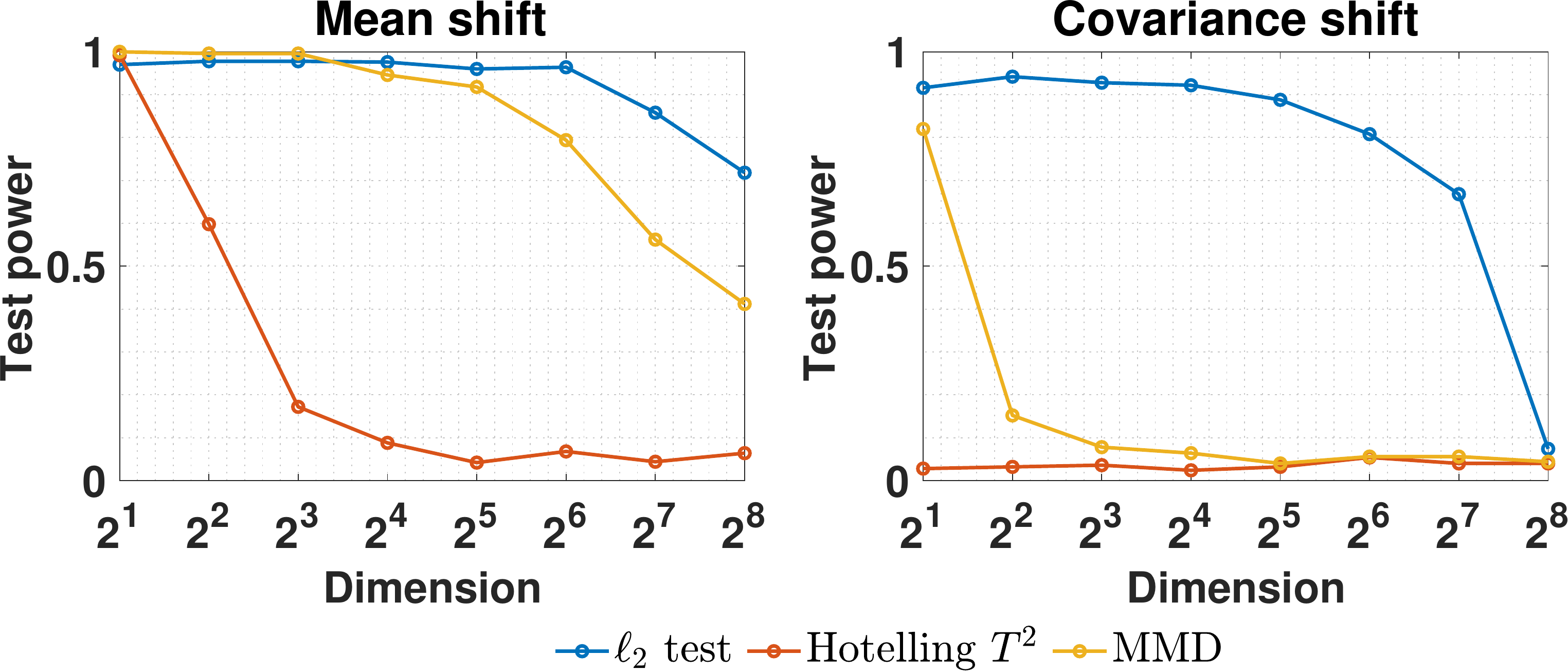}
\vspace{-0.1in}
\caption{Comparison of test power of the proposed test versus classic Hotelling's $T^2$ statistic and the MMD statistic, when performing a two-sample test on two Gaussian distributions, with significance level $\alpha=0.05$. (Left) Gaussian distributions having the same variance and different means; (Right) Gaussian distributions having same mean and different variances.}
\label{fig:two_sample_test_dimension}
\vspace{-0.3in}
\end{figure}

The test power for different methods is shown in Fig.~\ref{fig:two_sample_test_dimension}. The test power drops when the dimension increases, which is consistent with the results in \cite{ramdas2015decreasing}. Hotelling's $T^2$ test performs good in low dimensions, but its performance degrades quickly when we consider higher dimensional problems. The MMD test is comparable to $\ell_2$ test in low dimensions, but the $\ell_2$ test tends to outperform the MMD test in high dimensions. The reason can be that by projecting to one-dimensional spaces using a good projection, the power of $\ell_2$ test tends to decrease slower compared to Hotelling's $T^2$ and MMD tests. 

\subsection{Offline Change Detection}

As an extension and application of the proposed $\ell_2$ test, we investigate the performance for the {\it offline} change detection and compare the detection power, i.e., the probability of successfully detecting the change when there is a change. 

Assume we have sample $x_1,\ldots,x_{T}$ with a fixed time horizon $T=200$, when there is a change, we set the change-point $K=100$. The $\ell_2$ detection statistic at each time $t$ is $M\chi$ with $\chi$ defined in \eqref{chi-def} (here $M = 2LR/(L+R)$ is the normalizing constant). To avoid the segment being too short, we compute the detection statistics for time instances $t\in[w , T - w]$ with $w=20$, and then take the maximum. Similarly, the Hotelling's $T^2$ statistic at each time $t$ is computed using \eqref{eq:T2} by treating data before $t$ as one sample and after $t$ as another sample; the MMD statistic is computed in a similar way from \eqref{eq:MMD}. 
We claim there is a change-point when the maximum of the detection statistics within window $t\in[w, T - w]$ exceeds the threshold. The thresholds for different methods will be chosen by Monte Carlo simulation to control the false alarm rate.

We consider the following cases (distribution changes) in the numerical experiments.
\begin{enumerate}[label={\it Case \arabic*}, leftmargin=*,noitemsep]
\item (Discrete distributions). The support size is $n=10$, distribution shifts from $p = \bm{1}/10$ (uniform) to $q = [1/30,2/30,\ldots,5/30,5/30,\ldots,2/30,1/30]$ (non-uniform).
\item (Gaussian mean and covariance shift). The distribution shifts from two-dimensional standard Gaussian $\mathcal N(0,I_2)$ to $\mathcal N\left([0.5 \;\; 0]^T,[ 1 \;\;  0.7]^T [ 1 \;\; 0.7] +[ -1 \;\; 0.4 ]^T [ -1, \;\; 0.4 ]\right)$.
\item (Gaussian to Gaussian mixture). The distribution shifts from $\mathcal N(0, I_{20})$ to the Gaussian mixture $0.8\mathcal N(0, I_{20}) + 0.2\mathcal N(0, 0.1I_{20})$. 
\item (Gaussian to Laplace). The distribution shifts from standard Gaussian $\mathcal N(0,1)$ to Laplace distribution with zero mean and standard deviation $0.8$. 
\end{enumerate}

The detection power is averaged over 500 repetitions and is reported in Table \ref{tab:compare_power}. It shows that the proposed $\ell_2$ test outperforms the classic Hotelling's $T^2$ and MMD tests, especially when the distribution change is difficult to detect (such as Case 3 and Case 4, where pre- and post-change distributions are close). \textcolor{black}{For Case 2 detecting mean and covariance shifts, the MMD test performs slightly better. A possible explanation is that the MMD metric can capture the difference between pre- and post-change Gaussian distributions well in a fairly low-dimensional setting.}

\begin{table}[!h]\setlength\tabcolsep{3pt}
\renewcommand\arraystretch{1}
\centering
\caption{Detection power in offline change detection. The sequence of length is $200$. Thresholds for all methods are calibrated so that the significance level is $\alpha=0.10$ and $\alpha= 0.25$. Averaged over 500 trials.} \label{tab:compare_power}
\begin{tabular}{ccccc@{\hskip 0.25in}cccc}
\specialrule{.08em}{0em}{0em}
& \multicolumn{4}{c}{$\alpha=0.10$ } & \multicolumn{4}{c}{$\alpha = 0.25$} \\
& Case 1  & Case 2  & Case 3 & Case 4  & Case 1  & Case 2  & Case 3 & Case 4\\
\hline 
$\ell_2$ test 
& {\bf 0.52} & 0.85 & {\bf 0.18} &  {\bf 0.56}
& {\bf 0.70} & 0.90 & {\bf 0.35} & {\bf 0.71} \\
\hline 
MMD & 0.32 &{\bf 0.90} & 0.16 & 0.43
&0.60 & {\bf 0.95} & 0.34 & 0.69 \\
\hline
Hotelling's $T^2$ & 0.07 & 0.23 & 0.09 &  0.06
&  0.20 & 0.23 & 0.23  & 0.23 \\
\specialrule{.08em}{0em}{0em}
\end{tabular}
\vspace{-0.2in}
\end{table}

\subsection{Online Change Detection} \label{online_eg}

We further investigate the performance for {\it online} change detection and compare the average detection delay, i.e., the number of samples it takes to detect the change after the change happens. More specifically, the detection delay is the difference between the stopping time and the true change-point.

Assume we have samples $\{x_t,t=1,2,\ldots\}$ that are available sequentially. We adopt the convention that there are pre-change samples available as $\{x_0,x_{-1},\ldots\}$, which are referred as historical data and can be used during the detection procedure. Consider the window-limited $\ell_2$ detection procedure defined in \eqref{eq:stoptime_window} with parameter $m_0=20$ and $m_1=100$. The Hotelling's $T^2$ detection statistic at each time $t$ is constructed as 
$(\bar x_t - \hat\mu)^T\widehat\Sigma^{-1}(\bar x_t - \hat\mu)$
where $\bar x_t$ is the average of samples within window $[t-m_0+1,t]$, and $ \hat\mu,\widehat\Sigma$ are estimated from historical data. The MMD statistic is constructed in the same way as in \cite{li2015m} with block size $B_0=20$ and number of blocks $N=5$. We will claim change and stop the detection procedures when the detection statistic exceeds the threshold; the thresholds for different methods are chosen by Monte Carlo simulation to control the average run length.

\begin{figure}[h!]
\centering
\begin{tabular}{cccc}
\includegraphics[width = 0.23\textwidth]{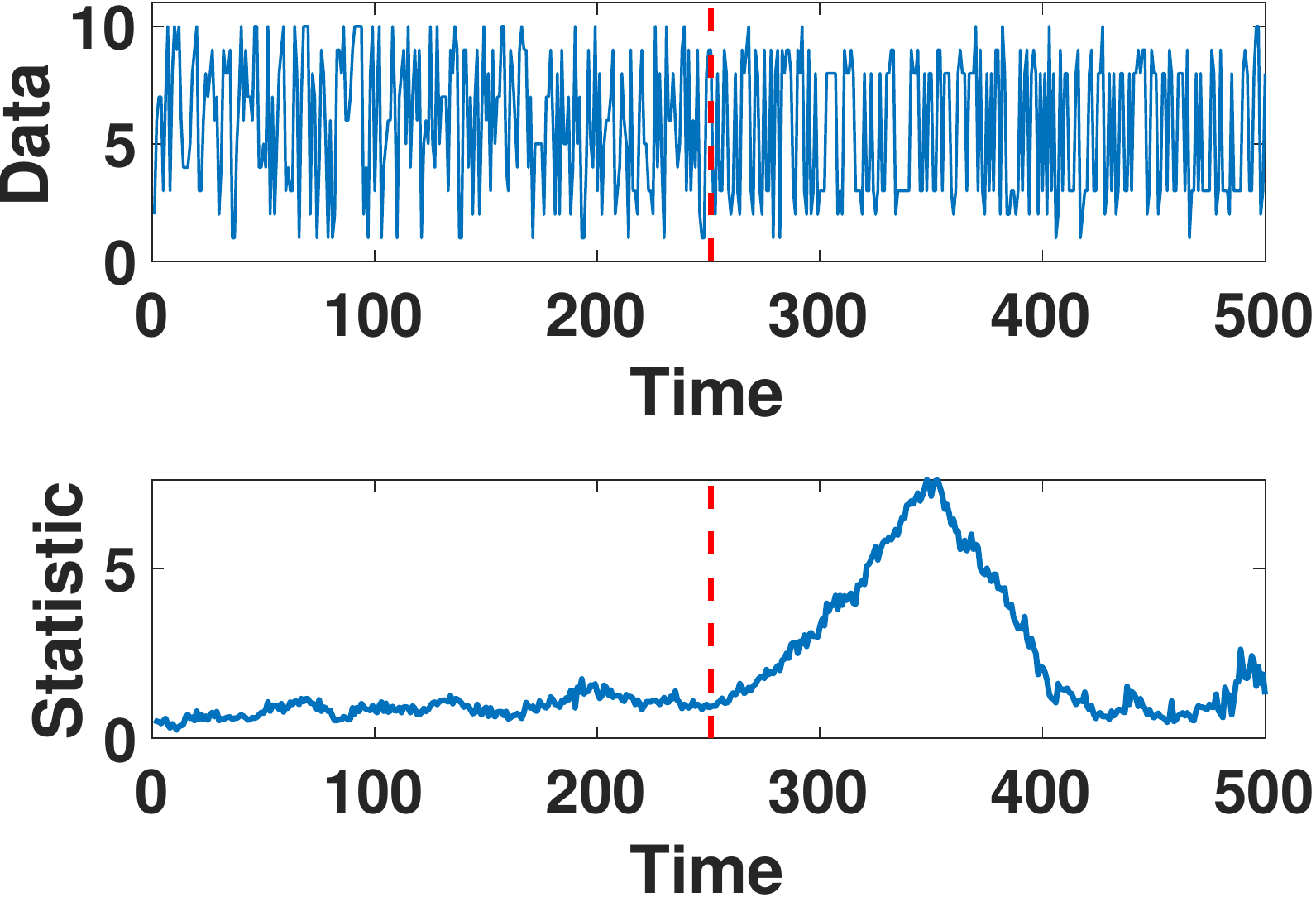} & \includegraphics[width = 0.23\textwidth]{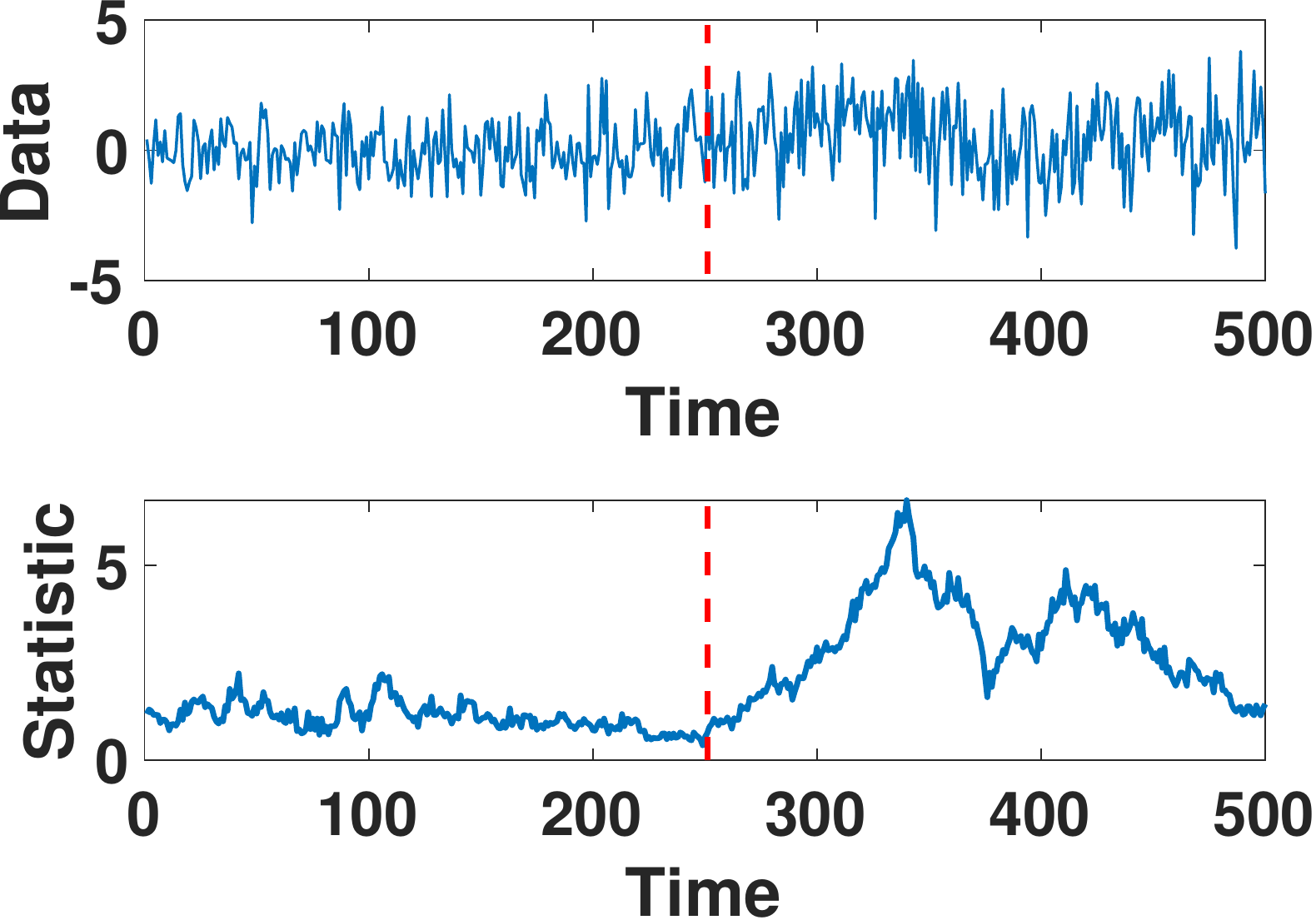} & \includegraphics[width = 0.23\textwidth]{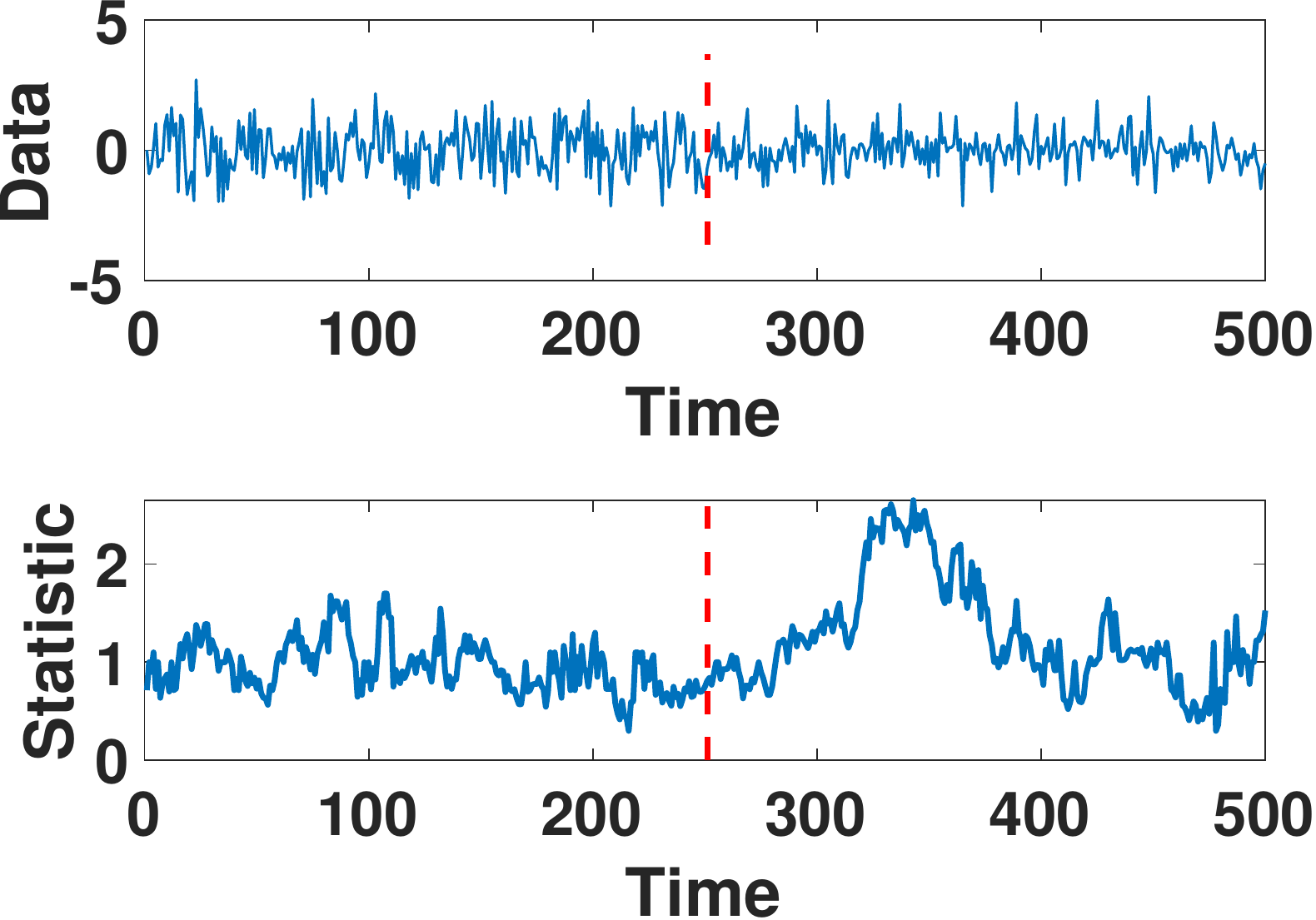}   &  \includegraphics[width = 0.23\textwidth]{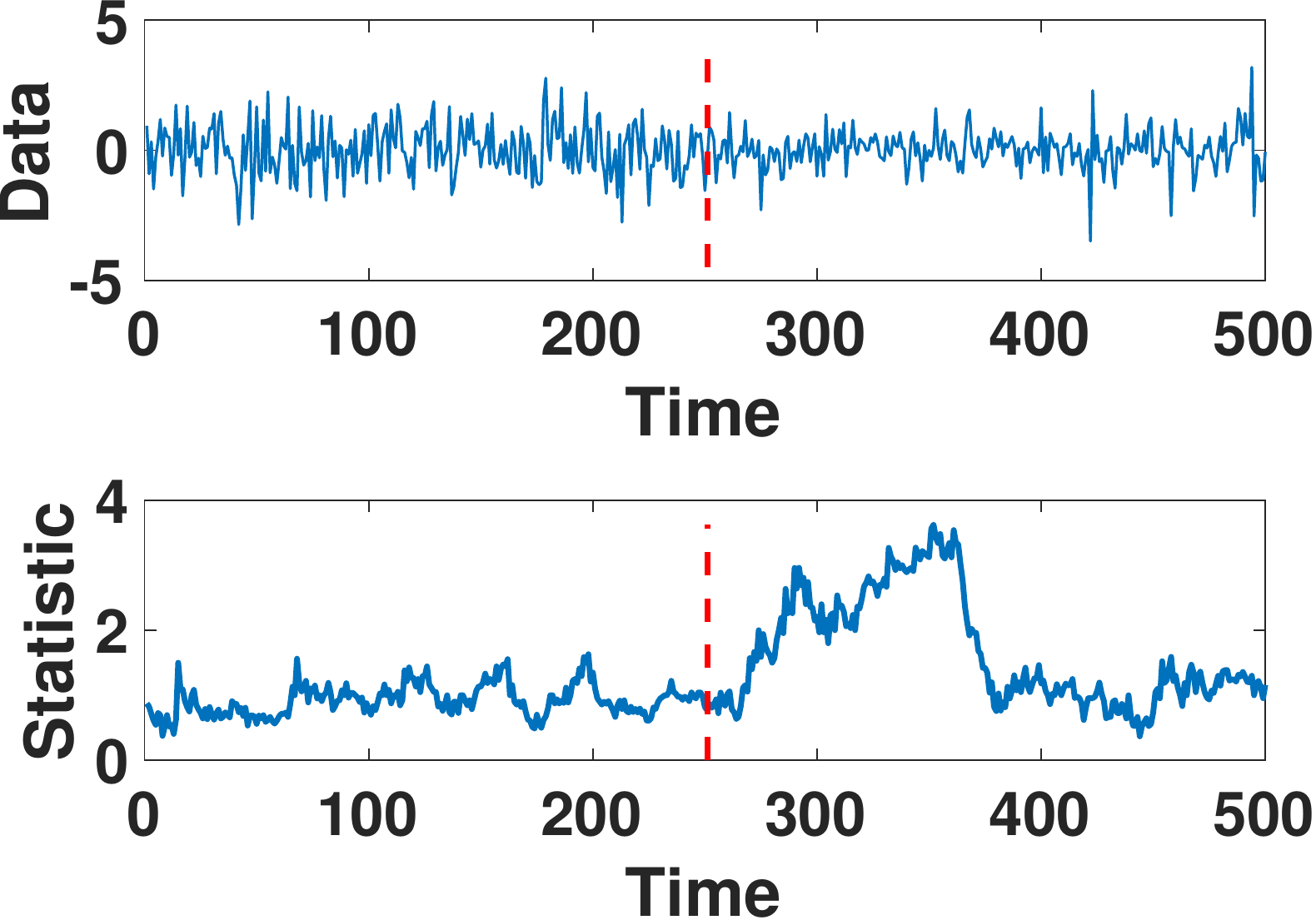} \\ 
{\small Case 1} & {\small Case 2} & {\small Case 3} & {\small Case 4} \\
\end{tabular}
\vspace{-0.1in}
\caption{Illustration of online change detection using the $\ell_2$ divergence under four simulated cases explained in Section \ref{online_eg}. For each case, the upper plot shows the raw data and the bottom plot shows the evolution path of the $\ell_2$ detection statistic, with true change-point indicated in red dash lines. }
\label{fig:detect_stat}
\vspace{-0.2in}
\end{figure}

We consider the following four cases, which are modified slightly from the offline case. We have increased the signal-to-noise ratio in certain cases to increase the detectability in the online setting.
\begin{enumerate}[label={\it Case \arabic*}, leftmargin=*,noitemsep]
\item (Discrete distributions). The support size is $n=10$, distribution shifts from $p = \bm{1}/10$ (uniform) to $q = [0.04,0.14,0.32,0,0,0,0,0.32,0.14,0.04]$ (non-uniform).  
\item (Gaussian mean and covariance shift). The distribution shifts from two-dimensional standard Gaussian $\mathcal N(0,I_2)$ to $\mathcal N\left([0.5 \;\; 0]^T,[ 1 \;\; 0.7 ]^T [1 \;\; 0.7]+[ -1 \;\; 0.4 ]^T [-1 \;\; 0.4 ] \right)$. 
\item (Gaussian to Gaussian mixture). The distribution shifts from $\mathcal N(0, I_{20})$ to the Gaussian mixture $0.4\mathcal N(0, I_{20}) + 0.6\mathcal N(0, 0.1I_{20})$. 
\item (Gaussian to Laplace). The distribution shifts from standard Gaussian $\mathcal N(0,1)$ to Laplace distribution with zero mean and standard deviation $0.7$. 
\end{enumerate}
The evolution paths of detection statistics for all cases are given in Fig.~\ref{fig:detect_stat}.
To simulate EDD, we let the change occur at the first time instant of the testing data. The detection delay is averaged over 500 repetitions and reported in Table\,\ref{tab:compare_EDD}.

\begin{table}[!h]\setlength\tabcolsep{12pt}
\renewcommand\arraystretch{1.0}
\vspace{-0.2in}
\centering
\caption{Comparison of EDD for online change detection using the proposed statistic, the MMD, and the Hotelling's $T^2$ statistic. The parameter is $n=10$, $m_0=20$, $m_1=100$ and thresholds for all methods are calibrated so that $\mbox{ARL}=500$. The dashed line indicates the method fails to detect the change (i.e., the delay is larger than the time horizon).}\label{tab:compare_EDD}
\begin{tabular}{ccccc}
\specialrule{.08em}{0em}{0em}
& Case 1  & Case 2  & Case 3 & Case 4 \\
\hline 
$\ell_2$ test & {\bf 20.34} & 89.66  & {\bf 69.23}  &  92.49 \\
\hline
MMD & 258.02 & 47.72 & --- & 394.91\\
\hline
Hotelling's $T^2$ & 406.42 & {\bf 36.79} &  --- &  370.61 \\
\specialrule{.08em}{0em}{0em}
\end{tabular}
\vspace{-0.2in}
\end{table}

\section{Real-data Study: Online gesture change detection}\label{sec:realdata}

In this section, we apply our method to the sequential gesture detection problem using a real dataset: the Microsoft Research Cambridge-12 (MSRC-12) Kinect gesture dataset \cite{fothergill2012instructing}. This dataset consists of sequences of human skeletal body part movements (represented as body part locations) collected from 30 people performing 12 gestures. There are 18 sensors in total, and each sensor records the coordinates in the three-dimensional Cartesian coordinate system at each time. Therefore there are 54 attributes, denoted by $y_t \in \bR^{54}$, $t=1,2,\ldots,T$. The goal is to detect the transition of gestures from the sequences of sensor observations. 

We apply the proposed online change detection procedure defined in \eqref{eq:stoptime_window} to the MSRC-12 dataset, and the detailed scheme is outlined as follows. We first preprocess the data by removing the frames that the person is standing still or with little movements. Then we select a {\it unit-norm} vector $u \in \bR^{54}$ and project data into this direction to obtain a univariate sequence: $x_t = u^\intercal y_t$. The projection vector $u$ is found by finding the optimal projection to maximize the Wasserstein distance described in Section~\ref{sec:optproj}. Then we discretize the univariate sequence into $n$ bins. At each time $t$, we construct the detection statistic $\max_{m_0 \leq t-k \leq m_1} \chi_{t,k}$ as illustrated in Fig.~\ref{fig:slidingwindow}. 

\begin{figure}[H]
\centering
\includegraphics[width = 0.45\textwidth]{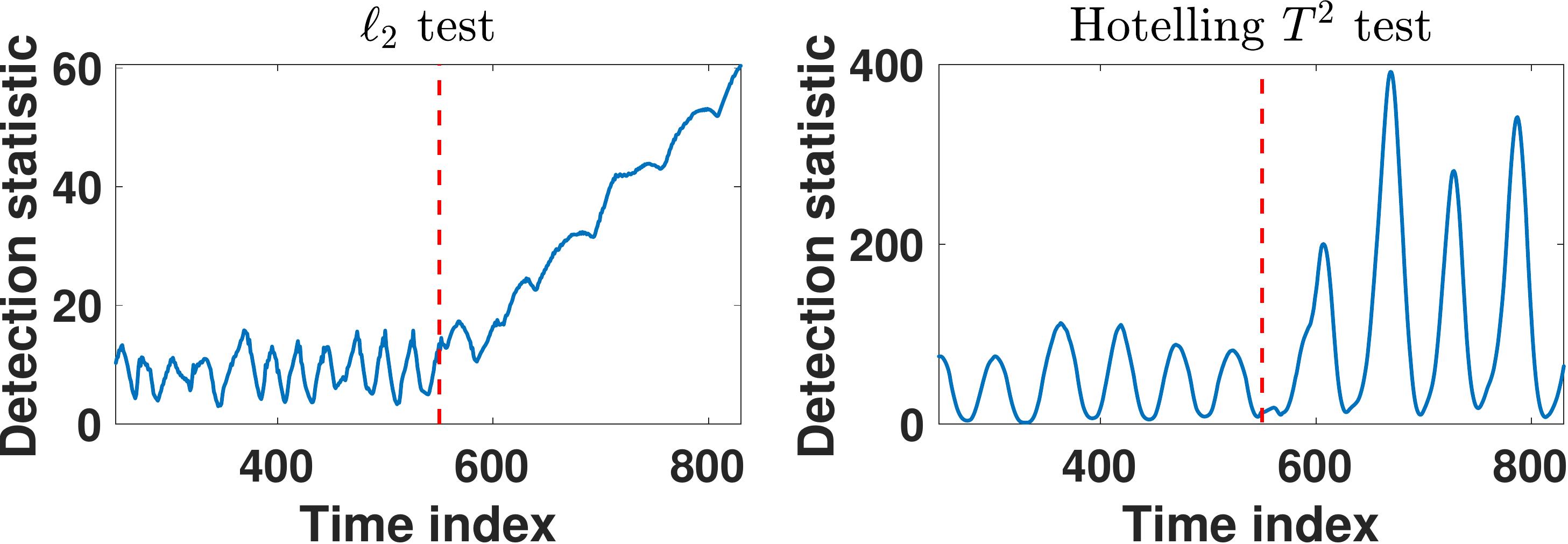}
\includegraphics[width = 0.45\textwidth]{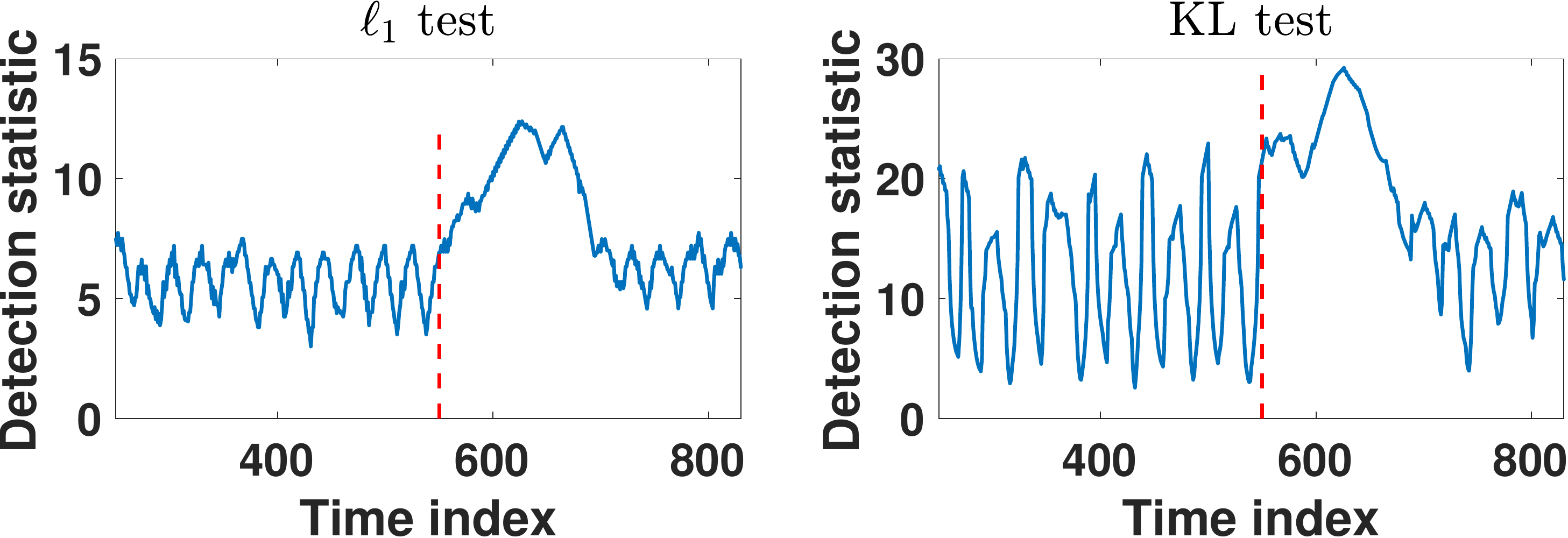}\\
\caption{Real-data example using online gesture change detection. Comparison of detection statistics (under uniform weights) for ``bow'' to ``throw'', for the proposed procedure, the Hotelling's $T^2$ test, $\ell_1$ test, and the KL test. Red dash lines indicate the true change-point (hand-labeled).}
\label{fig:random_optimal_stat}
\vspace{-0.2in}
\end{figure}

The parameters are set as $m_0 = 20$, $m_1 = 300$ for the detection procedure. The detection statistics are shown in Fig.~\ref{fig:random_optimal_stat}, with the true change indicated by red dash lines. 
We also compare the $\ell_2$ test with tests based on Hotelling's $T^2$ statistic, $\ell_1$ distance, and KL divergence. Using the $\ell_1$ approach detailed in Section\,\ref{sec:ell_pi}, we build the test statistic $\delta_{t,k}$ for time $t$ and potential change-point $k<t$, then the detection statistic is computed by maximizing over all potential change-points: $\max_{m_0 \leq t-k \leq m_1}\delta_{t,k}$. For the KL divergence test, at time $t$ and for $k<t$, set $N_{t,k} =\lceil t-k\rceil$ and denote the empirical distribution of two segments $[k-N_{t,k} +1 , k]$, $[k+1:t ]$ as $\hat f_1$ and $\hat f_2$ (with zero value adjusted to a small constant $\delta=0.01$), then compute their KL divergence as detection statistic.

%
% We can see that the performance of 
% %\subsubsection{Comparison of random projection and optimal projection}
% the random projection can be unsatisfactory and highly random. If we use the optimal projection method, two ``training'' samples are needed and we maximize the Wasserstein distance between projected observations. Under the sequential setting where we would like to make the decision in a timely manner, there is a trade-off between the estimation accuracy and the computational cost. 
%However, we observe that the estimation of projection $u$ does not differ much as $T$ varies. Therefore, we can obtain the optimal projection with fairly less ``training'' samples.
%Another trade-off exists for the choice of the length $L$ (assume $L=R$). The detection statistic $\chi$ has a large deviation when $L$ is small, while the detection delay will be large if $L$ is large, since we have to wait until enough post-change samples to raise the alarm.

From the results in Fig.~\ref{fig:random_optimal_stat}, we can observe that our sequential detection procedure based on the $\ell_2$ divergence can detect the change right after it happens. This is because the detection statistic before the change has a smaller variance, which indicates that we can set the threshold to be reasonably low for quicker detection. Moreover, there is a clear linear increasing trend after the change, enabling quick and reliable detection. In contrast, Hotelling's $T^2$ statistic does not have the desired online change detection behavior. The detection statistic is noisy before the change and does not have a consistent positive shift after the change; the KL test is even worse in this regard. The $\ell_1$ divergence-based test has a similar behavior as the $\ell_2$ divergence. However, the $\ell_1$ divergence has smaller ``signal-to-noise'' ratio in that the variance between the change is larger, and post-change distribution drift seems to be smaller.

\section{Conclusion}

We have presented a new non-parametric change detection procedure based on the optimal weighted $\ell_2$ divergence. We studied the optimality and various theoretical properties of the weighted $\ell_2$ divergence for the offline and online change-point detection. We also studied the practical aspects, including calibration threshold using training data, optimizing weights, and finding an optimal projection for high-dimensional data. We demonstrate the good performance of the proposed method using simulated and real-data for human gesture detection. 
%Note that we can have a situation where the empirical distributions live on different supports. 

%\textcolor{red}{Fundamentally, we are addressing the problem of testing the difference between two distributions based on samples only. Here, we present a simplified approach based on the $\ell_2$ divergence between empirical distributions and demonstrate its good performance compared with the commonly used distributional-free approach MMD test and the classical parametric approach Hotelling $T^2$ test.} 

% Acknowledgements should go at the end, before appendices and references

{\small \section*{Acknowledgement}

The authors are grateful to Professor Arkadi Nemirovski and Professor Anatoli Juditsky, and Professor Yihong Wu for their helpful discussions regarding optimality proof. We thank anonymous reviewers for their helpful and constructive comments. The authors would like to acknowledge support for this project from the National Science Foundation through NSF CAREER CCF-1650913, DMS-1938106, DMS-1830210, and CCF-1442635.}

%\clearpage
\bibliographystyle{IEEEtran}
\bibliography{ReferencesChP}

 \clearpage
\appendix
% \small
% \section{Proof of Lemma \ref{lem1002}}
\section{Proofs}\label{app}

\subsection{Proof of Proposition \ref{thm:concentration} (Test Properties)}\label{app:proof1}
Observe that the expectation of the empirical distribution of $N$-element sample drawn from a distribution $r\in\Delta_n$ on $\Omega$
is $r$, and the covariance matrix is \[
C_{r,N}={1\over N}[\Diag\{r\}-rr^T]\preceq{1\over N}\Diag\{r\}.
\]
Representing
$
\omega=p+\alpha,\omega'=p+\alpha',\zeta=q+\beta,\zeta'=q+\beta'$, we have
\[
\chi=(p-q)^T\Sigma(p-q)+ \underbrace{(\alpha-\beta)^T\Sigma(p-q)}_{B}+\underbrace{(\alpha'-\beta')^T\Sigma(p-q)}_{B'} + \underbrace{(\alpha-\beta)^T\Sigma(\alpha'-\beta')}_{C}.
\]
\noindent
Since $\alpha$ and $\beta$ are zero-mean and independent, the covariance matrix of $\alpha-\beta$ (and $\alpha'-\beta'$) is
\[
\begin{array}{c}
C_{p,L}+C_{q,R}=M^{-1}\left[2\gamma[\Diag\{p\}-pp^T]+2\bgamma[\Diag\{q\}-qq^T]\right]\preceq 2M^{-1}\Diag\{\gamma p+\bgamma q\},
\end{array}
\]
\noindent
whence
$$
\begin{array}{rcl}
\bE\{B^2\}&=&(p-q)^T\Sigma[C_{p,L}+C_{q,R}]\Sigma(p-q)\\
&\leq&(p-q)^T\Sigma[2M^{-1}\Diag\{\gamma p+\bgamma q\}]\Sigma(p-q)\\
&=&2M^{-1}\sum_i\sigma_i^2(p_i-q_i)^2(\gamma p_i+\bgamma q_i),\\
\end{array}
$$
and similarly $
\bE\{[B']^2\}\leq 2M^{-1}\sum_i\sigma_i^2(p_i-q_i)^2(\gamma p_i+\bgamma q_i)$. 
Moreover,
\[
\begin{array}{rcl}
\bE\{C^2\}&=&\bE\left\{\left[\sum_i\sigma_i(\alpha_i-\beta_i)(\alpha^\prime_i-\beta^\prime_i)\right]^2\right\}\\
% &=\sum_{i,j}\bE\{\sigma_i\sigma_j(\alpha_i-\beta_i)(\alpha^\prime_i-\beta^\prime_i)(\alpha_j-\beta_j)(\alpha^\prime_j-\beta^\prime_j)\}\\
&=&\sum_{i,j}\sigma_i\sigma_j
\bE\{(\alpha_i-\beta_i)(\alpha_j-\beta_j)\}\bE\{(\alpha^\prime_i-\beta^\prime_i)(\alpha^\prime_j-\beta^\prime_j)\}\\
% & \quad\quad\quad\hbox{\ [since $^\prime$-variables are independent of non-$^\prime$ ones]}\\
&=&\sum_{i,j}\sigma_i\sigma_j [C_{p,L}+C_{q,R}]_{ij}^2\\
&=&\sum_i\sigma_i^2[C_{p,L}+C_{q,R}]_{ii}^2+\sum_{i\neq j}\sigma_i\sigma_j[C_{p,L}+C_{q,R}]_{ij}^2\\
&=&4M^{-2}\left[\sum_i\sigma_i^2\left[\gamma [p_i-p_i^2]+\bgamma[q_i-q_i^2]\right]^2+\sum_{i\neq j}\sigma_i\sigma_j
[\gamma p_ip_j+\bgamma q_iq_j]^2\right]\\
&\leq&
4M^{-2}\left[\sum_i\sigma_i^2[\gamma p_i+\bgamma q_i]^2+\sum_{i\neq j}\sigma_i\sigma_j[\gamma p_ip_j+\bgamma q_iq_j]^2\right]\\
&\leq& 4M^{-2}\left[\gamma\sum_i\sigma_i^2p_i^2+\bgamma \sum_i\sigma_i^2q_i^2+\gamma \sum_{i\neq j}\sigma_i\sigma_jp_i^2p_j^2 +\bgamma \sum_{i\neq j}\sigma_i\sigma_jq_i^2q_j^2\right]\\
&\leq& 4M^{-2}\left[\gamma\sum_i\sigma_i^2p_i^2+\bgamma\sum_i\sigma_i^2q_i^2+
\gamma \left[\sum_i\sigma_ip_i^2\right]^2+ \bgamma \left[\sum_i\sigma_iq_i^2\right]^2\right].
\end{array}
\]
\noindent
Note that
$$
\begin{array}{c}
[\sum_i\sigma_ip_i^2]^2=[\sum_i(\sigma_ip_i)p_i]^2\leq [\sum_i\sigma_i^2p_i^2][\sum_ip_i^2]\leq \sum_i\sigma_i^2p_i^2,
\end{array}
$$
which combines with the previous computation to imply that
$$
\begin{array}{c}
\bE\{C^2\}\leq 8M^{-2}\left[\gamma\sum_i\sigma_i^2p_i^2+\bgamma\sum_i\sigma_i^2q_i^2\right].
\end{array}
$$
Consequently, by applying Chebyshev's inequality to $B$, $B'$, and $C$, respectively, for every $\theta\geq1$, the probability $\pi(\theta)$ of the event
$$
\begin{array}{rl}
&|\chi-\sum_i\sigma_i(p_i-q_i)^2| \leq  \\
&\hspace{80pt} 2\sqrt{2}\theta[\sqrt{M^{-1}\sum_i\sigma_i^2(p_i-q_i)^2(\gamma p_i+\bgamma q_i)} + M^{-1}\sqrt{\gamma\sum_i\sigma_i^2p_i^2+\bgamma\sum_i\sigma_i^2q_i^2}]
\end{array}
$$
% {%\small
% \begin{multline*}
% |\chi-\sum_i\sigma_i(p_i-q_i)^2| \leq \\
% 2\sqrt{2}\theta\left[\sqrt{M^{-1}\sum_i\sigma_i^2(p_i-q_i)^2(\gamma p_i+\bgamma q_i)} + M^{-1}\sqrt{\gamma\sum_i\sigma_i^2p_i^2+\bgamma\sum_i\sigma_i^2q_i^2}\right]
% \end{multline*}}
\noindent
is at least $1-3/\theta^2$. By inspecting the derivation, it is immediately seen that when $p=q$, the lower bound $1-3/\theta^2$ on $\pi(\theta)$ can be improved to $1-1/\theta^2$, so that
\begin{equation*}\label{pigeq}
\pi(\theta)\geq\left\{\begin{array}{ll}1-3/\theta^2,&p\neq q,\\
1-1/\theta^2,&p=q.\\
\end{array}\right.
\end{equation*}
Taking into account what $\cT$ is, the conclusion of Proposition~\ref{thm:concentration} follows.
\begin{remark}
Note we do not use the standard ``Poissonization'' approach which assumes that, rather than drawing $N$ independent
samples from a distribution, first select $N'$ from Poisson distribution with mean value $N$, and then draw $N'$ samples. Such Poissonization makes the number of times different elements occur in the sample independent, simplifying the analysis. Instead, we model the empirical distribution directly by considering the dependence in the covariance matrix $C_{r,N}$. 
\end{remark}
\subsection{Proof of Proposition \ref{prop:estimate} (Bounding $\|p\|_2^2$)}

\noindent
\textbf{Step 1.} Let $\xi$, $\xi'$ be empirical distributions of observations in two consecutive $Q$-element segments of sample $X^1$ that are generated from distribution $p$. 
Setting $\eta=\xi-p$, $\eta'=\xi'-p$, we have
$$
\xi^T\xi'=p^Tp+\underbrace{p^T\eta}_{B}+\underbrace{p^T\eta'}_{B'}+\underbrace{\eta^T\eta'}_{C}.
$$
Since $\eta$ and $\eta'$ are zero-mean vectors and independent of each other, with covariance matrix $C_{p,Q}=Q^{-1}[\Diag\{p\}-pp^T]\preceq Q^{-1}\Diag\{p\}$, we have
$$
\begin{array}{c}
\bE\{B^2\}=\bE\{[B']^2\}=p^T C_{p,Q}p\preceq Q^{-1}p^T\Diag\{p\}p\leq Q^{-1}\sum_ip_i^3,
\end{array}
$$
and
$$
\begin{array}{rcl}
\bE\{C^2\}&=&\sum_{i,j}\bE\{\eta_i\eta_i^\prime\eta_j\eta_j^\prime\}\\
&=&\sum_{i,j}\left[\bE\{\eta_i\eta_j\}\right]^2
=Q^{-2}\left[\sum_i[p_i-p_i^2]^2+\sum_{i\neq j}p_i^2p_j^2\right]
\leq2Q^{-2}\sum_ip_i^2.
\end{array}
$$
Consequently, by applying Chebyshev's inequality to $B$, $B'$, and $C$, respectively, the probability of the event
\begin{equation}\label{pp}
\begin{array}{c}
|\xi^T\xi'-p^Tp|>3\left[2Q^{-1/2}\sqrt{\sum_ip_i^3}+\sqrt{2}Q^{-1}\|p\|_2\right]
\end{array}
\end{equation}
is $\leq 1/3$. Taking into account that $\sum_ip_i^3\leq \|p\|_2^2\|p\|_\infty\leq \|p\|_2^3$, we have proved the following

\begin{lemma}[Concentration Inequality for $\xi^T\xi'$]\label{lemnew} Assume that there exists $\rho \in \mathbb R_{+}$ such that the distribution $p$ satisfies the relation
\begin{equation*}\label{letrho}
\|p\|_2\leq \sqrt{2}\rho,
\end{equation*}
and a positive integer $Q$ be such that
% \begin{equation}\label{suchthat}
% 3\left[2^{5/2}Q^{-1/2}\rho^3+2Q^{-1}\rho\right]\leq {1\over 3}\rho^2.
% \end{equation}
\begin{equation}\label{suchthat}
3\left[2^{7/4}Q^{-1/2}\rho^{3/2}+2Q^{-1}\rho\right]\leq {1\over 3}\rho^2.
\end{equation}
When $\xi$, $\xi'$ are empirical distributions of two consecutive segments, of cardinality $Q$ each, generated from distribution $p$, we have
\begin{equation}\label{emma}
\bP\left\{|\xi^T\xi'-\|p\|_2^2|>{1\over 3}\rho^2\right\}\leq {1 \over 3}.
\end{equation}
\end{lemma}

\noindent
\textbf{Step 2.} The parameters $Q=Q_i$ and $\rho=\rho_i$ of $i$-th stage of the training-step in Algorithm~\ref{alg:pilot} satisfy (\ref{suchthat}). Recalling that by the definition of $i(p)$ we have $\|p\|_2\leq \sqrt{2}\rho_i$, $1\leq i\leq i(p)$. Invoking Lemma \ref{lemnew} and the definition of $S$ in \eqref{eq:S}, we conclude that the probability of the event
$$
\cE: |\Theta_i-\|p\|_2^2|\leq {1\over 3}\rho_i^2,\,\forall 1\leq i\leq i(p),
$$
is at least $1-\delta$. Assume that this event takes place. By the definition of $i(p)$, we have $\rho^2_{i(p)}\leq\|p\|_2^2$, and since we are in the case of $\cE$, we have also
$|\Theta_{i(p)}-\|p\|_2^2|\leq {1\over 3}\rho_{i(p)}^2$, whence
$$
\Theta_{i(p)}\geq {2\over 3}\|p\|_2^2\geq {2\over 3}\rho_{i(p)}^2.
$$
We see that under our assumption the trial run ends up with a success at some stage $k\leq i(p)$, so that
$$
\Theta_k\geq {2\over 3}\rho_k^2\ \text{and} \ |\Theta_k-\|p\|_2^2|\leq {1\over 3}\rho_k^2
$$
(the second relation holds true since we are in the case of $\cE$). As a result,
$$
\varrho=\Theta_k+{1\over 3}\rho_k^2\geq \|p\|_2^2,
$$
and $|\Theta_k-\|p\|^2_2|\leq \rho_k^2/3 \leq \Theta_k/2$, implying that $\Theta_k\leq 2\|p\|^2$, whence
$$
\varrho=\Theta_k+{1\over 3}\rho_k^2\leq \Theta_k + \frac12 \Theta_k = {3\over 2}\Theta_k\leq 3\|p\|_2^2.
$$
We see that, with probability at least $1-\delta$, the number of stages in the training-step is at most $i(p)$, and the output $\varrho$ of the test satisfies
(\ref{target}). Besides this, from \eqref{suchthat} it is immediately seen that $Q_i\leq O(1)/\rho_i$, so that $Q_{i(p)}\leq O(1)/\|p\|_2$ due to $\rho_{i(p)}\geq \|p\|_2/\sqrt{2}$. Thus, when the training-step stops before or at stage $i(p)$, the total number of observations used in training-step
indeed does not exceed $4SQ_{i(p)}\leq O(1)S/\|p\|_2$. Note that by the definition of $S$ in \eqref{eq:S}, we have
\begin{equation}\label{Sis}
S\leq O(1)\ln(\ln(n)/\delta),
\end{equation}
which implies (\ref{nini}).

\subsection{Proof of Theorem \ref{atmostbeta} (Sample Complexity)}

By Proposition \ref{prop:estimate}, with properly selected $O(1)$ in (\ref{condI}), the probability for the training-step to be successful is at least $1-\delta$, and there is enough observations to perform the $K$ individual tests of the testing stage. From now on we assume that $O(1)$ in (\ref{condI}) meets this requirement.

For $k\leq i(p)$, let $\cE_k$ be the condition stating that the training-step is successful and terminates at stage $k$. Note that this is a condition on the first $N_k=4SQ_k$ observations of the sample set $X^1$. Let us fix a realization of these $N_k$ observations satisfying condition $\cE_k$; from now on, speaking about probabilities of various events, we mean probabilities taken with respect to conditional, the above realization given, probability distribution of the remaining $N-N_k$ observations in sample $X^1$ and the entire $N$ observations in $X^2$. 

We first prove the type-I risk is at most $\alpha$. Note that we are in the situation when the training-step was successful, hence $\varrho \geq\|p\|_2^2$. Consequently,
the threshold (\ref{resultingell}) satisfies relation (\ref{ell}) with $\theta=3$ and $\sigma_i\equiv 1$, implying by the first claim in Proposition \ref{thm:concentration} (where we set $L=R=M$) that when $p=q$, the probability to claim $H_1$ by a particular one of the $K$ individual tests is at most $1/9$. By the definition of $K(\alpha,\beta)$, we conclude that the type-I risk is indeed at most $\alpha$.

We then prove the type-II risk is at most $\beta$ whenever the condition (\ref{Mislarge}) holds. Assume that $\|p-q\|_2\geq\epsilon\|p\|_2$ with some $\epsilon>0$, and set $L=R=M$, $\theta=3$ and $\sigma_i\equiv1$. With $\ell$ given by (\ref{resultingell}), the inequality (\ref{magnitude}) reads
\begin{equation}\label{magnitudeagain}
\begin{array}{rcl}
\|p-q\|_2^2&>& 6\sqrt{2}\big[M^{-1}\sqrt{\varrho}+M^{-1/2}\sqrt{\sum_i(p_i-q_i)^2(p_i+q_i)/2}+\\
&&\hspace{30pt}M^{-1}\sqrt{(\|p\|_2^2+\|q\|_2^2)/2}\big].
\end{array}
\end{equation}
Note that the condition (\ref{magnitudeagain}) ensures that the power of every individual test is at least $2/3$; thus, due to the choice of $K(\alpha,\beta)$, the type-II risk is at most $\beta$. It only remains to verify that condition (\ref{Mislarge}) implies the
validity of (\ref{magnitudeagain}). Since we are in the situation that the training-step is successful, the condition (\ref{target}) holds and in particular, $\|p\|_2\leq\sqrt{\varrho}\leq \sqrt{3}\|p\|_2$, implying that the right hand side in (\ref{magnitudeagain}) is at most
\[
\begin{array}{c}
{\cal R}=O(1)\left[M^{-1}[\sqrt{\varrho}+\|q\|_2]+M^{-1/2}\|p-q\|_2\sqrt{\sqrt{\varrho}+\|q\|_2}\right],
\end{array}
\]
and therefore in order to ensure the validity of (\ref{magnitudeagain}), it suffices to ensure that
\begin{equation}\label{toensure}
\|p-q\|_2^2\geq O(1)M^{-1}[\sqrt{\varrho}+\|q\|_2].
\end{equation}
First consider the case when $\|q\|_2\leq 2\|p\|_2$, which combines with (\ref{target}) to imply that the right hand side in (\ref{toensure})  is $\leq O(1)M^{-1}\sqrt{\varrho}$. By (\ref{target}) and $\|p-q\|_2\geq\epsilon\|p\|_2$, the left hand side in (\ref{toensure}) is at least $O(1)\epsilon^2\varrho$, so that (\ref{toensure}) indeed is implied by (\ref{Mislarge}), provided that the absolute constant factor in the latter relation is selected properly. Then consider the case when $\|q\|_2\geq 2\|p\|_2$. In this case, by (\ref{target}),
the right hand side in (\ref{toensure}) is at most $O(1)M^{-1}\|q\|_2$, and the left hand side in (\ref{toensure}) is at least $O(1)\|q\|_2^2$, implying that (\ref{toensure}) holds true when $M\geq O(1)/\|q\|_2$. The validity of the latter condition, in view of $\|q\|_2\geq2\|p\|_2$, clearly is guaranteed by the validity of (\ref{Mislarge}), provided $O(1)$ in the latter relation is selected properly.
 \qed

\subsection{Proof of Proposition \ref{thm:optimality} (Optimality)}
Assuming $n$ is even, consider the following two scenarios on distributions $p,q$ from which the two sets of sample $X^1$ and $X^2$ are independently generated:
\begin{enumerate}
\item Both samples are i.i.d. drawn from the uniform distribution on $\Omega$.
\item The nature draws, independently of each other, two $n/2$-element subsets, $\Omega_1$ and $\Omega_2$, of $\Omega$,
from the uniform distribution on the family of all subsets of $\Omega$ of cardinality $n/2$; the $N$-observation samples $X^1$ are i.i.d. drawn from the uniform distribution on $\Omega_1$, and the $N$-observation samples $X^2$ are i.i.d. drawn from the uniform distribution on $\Omega_2$.
\end{enumerate}
In the first scenario, the hypothesis $H_0$ is true; in the second, there is a {\it significant} difference between $p$ and $q$ -- with probability close to 1 when $n$ is large enough, we have the $\|p-q\|_2\geq \epsilon \|p\|_2 = \epsilon\sqrt{2/n}$ for any $\epsilon$ small enough, e.g., for $0<\epsilon<1/2$. Denote the union of two sets of samples as $x^{2N}:=(x_1,...,x_{2N})$. It follows that if there exists a test $\cT$ obeying the premise of Proposition~\ref{thm:concentration}, then there exists a low risk test deciding on whether the entire $2N$-element sample $x^{2N}$ shown to us is generated according to the first or the second scenario. 

Specifically, given $x^{2N}$, let us split it into two halves and apply to the two resulting $N$-observation samples the test $\cT$; if the test claim $H_1$, we conclude that $x^{2N}$ is generated according to the second scenario, otherwise we claim that $x^{2N}$ is generated according to the first scenario. When $x^{2N}$ is generated by the first scenario, the probability for $\cT$ to claim $H_1$ is at most $0.1$, that is, the probability to reject the first scenario when it is true is at most $\alpha$. On the other hand, when $x^{2N}$ is generated according to the second scenario, the conditional, $p$ and $q$ given, probability for $\cT$ to accept $H_0$ should be at most $\beta$, provided that $\|p-q\|_2\geq \epsilon\|p\|_2$ for a given $\epsilon\in(0,1/2)$; when $n$ is large, the probability for the condition $\|p-q\|_2\geq \epsilon\|p\|_2$ to hold true approaches 1, so that for large enough values of $n$, the probability for the condition to hold is at least $1-\beta$ and therefore the probability of claiming a sample generated according to the second scenario as one generated according to the first one, is at most $2\beta$. Thus, for properly selected $n_0$ and all $n\geq n_0$, given $x^{2N}$, we can decide with risk $\leq 2\beta$ on the scenario resulted in $x^{2N}$.

On the other hand, consider the distribution of $x^{2N}$. The corresponding observation space is the space $\Omega^{2N}$
of $2N$-element sequences with entries from $\Omega$. Let $\widehat{\Omega}$ be the part of $\Omega^{2N}$ comprised of sequences with {\it all entries different from each
other}, and $\widetilde{\Omega}$ be the complement of $\widehat{\Omega}$ in $\Omega^{2N}$. Let also $P_1$ and $P_2$ be the distributions of our observations under the
first and the second scenarios, and $P$ be the distribution on $\Omega^{2N}$ which assigns equal masses to all points from $\widehat{\Omega}$ and zero masses to the points
outside of $\widehat{\Omega}$. By evident symmetry reasons, we have
$$
P_i=(1-\epsilon_i)P+\epsilon_i Q_i,\,i=1,2,
$$
where $Q_1$ and $Q_2$ are probability distributions supported on $\widetilde{\Omega}$, and $\epsilon_i$ is the probability, under scenario $i$, to observe $2N$-element sample in $\widetilde{\Omega}$. We clearly have
$$
\epsilon_1\leq {N(2N-1)\over n},\,\epsilon_2\leq {4N(2N-1)\over n}.
$$
Indeed, for a fixed pair of indexes $t_1,t_2$, $1\leq t_1<t_2\leq 2N$, the probability to get $x_{t_1}=x_{t_2}$ in $x^{2N}$ is $1/n$ under
the first scenario and is at most $4/n$ under the second scenario, while the number of pairs $t_1,t_2$ in question is $N(2N-1)$. We see that
$$
\begin{array}{l}
\sum_{\zeta^{2N}\in\Omega^{2N}}\min\{P_1(\zeta^{2N}),P_2(\zeta^{2N})\}
=\min\{1-\epsilon_1,1-\epsilon_2\}\sum_{\zeta^{2N}\in\widehat{\Omega}}P(\zeta^{2N})\\
\multicolumn{1}{r}{+
\sum_{\zeta^{2N}\in\widetilde{\Omega}}\min\{\epsilon_1 Q_1(\zeta^{2N}),\epsilon_2Q_2(\zeta^{2N})\}\geq \min\{1-\epsilon_1,1-\epsilon_2\},}\\
\end{array}
$$
implying that our scenarios cannot be decided upon with risk $\leq\min\{1-\epsilon_1,1-\epsilon_2\}/2$. When $N\ll \sqrt{n}$, both $\epsilon_1$ and $\epsilon_2$ are small, so that it is impossible to decide on our scenarios with risks $\alpha$ (and $2\beta$). The bottom line is that under the premise of Proposition~\ref{thm:concentration}, we either have $n\leq n_0$, or $N\geq O(1)\sqrt{n}$ with properly selected constant $O(1)$ that depends on $\alpha,\beta,\epsilon$, and the conclusion follows. Note that here the dependence on reliability parameters is logarithmic and is not our focus here since here we aim to study the sample optimality with respect to the cardinality $n$ of the set $\Omega$. \qed

\subsection{Proof of Theorem~\ref{ell1_risk} (Properties of Generalized $\pi$-test)}

\noindent
\textbf{Step 1.}
We can upper-bound the Type-I risk as follows. Let $p$, $q$ be the underlying distributions from which $X^1$ and $X^2$ are sampled. Assume that $\tau=\tau_{N}'(\alpha)$ solves (\ref{equation}). The observations used to build $\hat{f}[x^\ts]$ and $\hat{f}[\bar{x}^\ts]$ are {\it independent} of those used to build $h$, so that the conditional, given $h$, probability $\epsilon$ of the event $h^T\hat{f}[x^\ts]>h^Tp+\tau$ is at most $\alpha/2$. Indeed, let  $p_i$ and $h_i$ be $i$-th entries of $p$ and $h$, respectively. Let $\xi$ be a random variable taking values $h_i-h^Tp$ with probability $p_i$ for $i=1,\ldots,n$, and let $\xi^1,...,\xi^N$ be  independent copies of $\xi$. The distribution of the random variable $h^T\hat{f}[x^\ts]-h^Tp$ is exactly the same as the distribution of the random variable ${1\over N}\sum_{t=1}^N\xi^t$. Taking into account that $\xi^t$ are i.i.d. zero mean random variables with variance $\leq \sum_ip_ih_i^2\leq 1$ taking values in $[-2,2]$ (recall that $\|h\|_\infty=1$) and applying Bernstein inequality, for every $s>0$ we have
$$
\begin{array}{c}
\bP\{h^T\hat{f}[x^\ts]-h^Tp>s\}=\bP\left\{\sum_{t=1}^N\xi^t>Ns\right\}\leq \exp\left\{-{Ns^2\over 2+{4\over 3}s}\right\}.
\end{array}
$$
Setting $s=\tau$, we have the above probability is upper bounded by $\alpha/2$.

By the same argument, the conditional, given $h$, probability of the event $h^T\hat{f}[\bar{x}^\ts]<h^Tq-\tau$ is $\leq\alpha/2$, implying that when $p=q$, the conditional, $h$ given, probability of the event $h^T[\hat{f}[x^\ts]-\hat{f}[\bar{x}^\ts]]>\ell:=2\tau$ is at most $\alpha$. In other words, we have proved that the type-I risk of our test is $\leq\alpha$, provided that $\tau$ solves (\ref{equation}). The only other option is that $\tau=1$, in which case the type-I risk is zero due to $\|h\|_\infty=1$ and $\|\hat{f}[x^\ts]-\hat{f}[\bar{x}^\ts]\|_1\leq 2$.
%\par
%Note that the fact that with threshold (\ref{ellthresh}) the Type-I risk of our test is $\leq\alpha$ is independent of what exactly $h$ is; all what matters is that $h$ is deterministic function of training sample such that $\|h\|_\infty\leq1$.

\noindent
\textbf{Step 2.}  Consider the ``bad training events''
\def\cB{{\cal B}}
\def\cG{{\cal G}}
$$
\begin{array}{rcl}
\cB^{\tr,1}&=&\{x^\tr:\|\hat{f}[x^\tr]-p\|_\pi>\rho_{\beta,n,\pi}(N)\};\\
\cB^{\tr,2}&=&\{\bar{x}^\tr:\|\hat{f}[\bar{x}^\tr]-q\|_\pi>\rho_{\beta,n,\pi}(N)\}.
\end{array}
$$
Note that all these events are described solely in terms of full training sample $(x^\tr,\bar{x}^\tr)$ and have probabilities $\leq\beta$ each. Let the ``good training event" $\cG^{\tr}$ be the set of all realizations of the full training sample where no one of the bad training events takes place, so that the probability of the good training event is at least $1-2\beta$.
Assume that a realization of the full training sample belongs to $\cG^{\tr}$, implying that
\begin{equation}\label{eq18}
\left\|\left[\hat{f}[x^\tr]-\hat{f}[\bar{x}^\tr]\right]-[p-q]\right\|_\pi= \left\|\left[\hat{f}[x^\tr]-p\right]-\left[\hat{f}[\bar{x}^\tr]-q\right]\right\|_\pi \leq 2 \rho_{\beta,n,\pi}(N).
\end{equation}
\noindent
Consequently,
\begin{equation}\label{eqA2}
\|\hat{f}[x^\tr]-\hat{f}[\bar{x}^\tr]\|_\pi\geq \|p-q\|_\pi-2\rho_{\beta,n,\pi}(N).
\end{equation}
Besides this,
denoting by $e$ the separator built at ``training stage'' and setting $\theta=\|e\|_\infty$, so that $\theta\leq1$ and $h=e/\theta$, we have
$$
\begin{array}{rcl}
h^T[p-q] &=& h^T\left[\hat{f}[x^\tr]-\hat{f}[\bar{x}^\tr]\right] + h^T\left[[p-q]-\left[\hat{f}[x^\tr]-\hat{f}[\bar{x}^\tr]\right]\right] \\
&\geq& h^T\left[\hat{f}[x^\tr]-\hat{f}[\bar{x}^\tr]\right]-\|h\|_{\pi_*}
\left\|[p-q]-\left[\hat{f}[x^\tr]-\hat{f}[\bar{x}^\tr]\right]\right\|_\pi\\
&\geq&\theta^{-1}\left[\|\hat{f}[x^\tr]-\hat{f}[\bar{x}^\tr]\|_\pi-2\rho_{\beta,n,\pi}(N)\right],
\end{array}
$$
which combines with (\ref{eqA2}) to imply that whenever the full training sample belongs to $\cG^\tr$, it holds
\begin{equation}\label{eq21}
h^T[p-q]\geq\theta^{-1}\left[\|p-q\|_\pi-4\rho_{\beta,n,\pi}(N)\right].
\end{equation}
The same reasoning applied to testing samples in the role of training ones says that given a realization of the full training sample, one can specify a set $\cG^{\ts}$ of realizations of the full testing sample such that the conditional, training sample given, probability for testing sample to belong to $\cG^\ts$ is $\geq1-2\beta$ and whenever this inclusion takes place, one has
\begin{equation}
\|[\hat{f}[x^\ts]-\hat{f}[\bar{x}^\ts]]-[p-q]\|_\pi\leq 2\rho_{\beta,n,\pi}(N).
\end{equation}
It follows that when the ``good" event $\cG$ --- the full training sample belongs to $\cG^\tr$, and the full testing sample belongs to $\cG^{\ts}$ --- takes place, we have
$$
h^T[\hat{f}[x^\ts]-\hat{f}[\bar{x}^\ts]]\geq h^T[p-q]-2\theta^{-1}\rho_{\beta,n,\pi}(N).
$$
We conclude that when good event takes place (which happens with probability at least $1-4\beta$), it holds
$$
h^T[\hat{f}[x^\ts]-\hat{f}[\bar{x}^\ts]]\geq \theta^{-1}\left[\|p-q\|_\pi-6\rho_{\beta,n,\pi}(N)\right],
$$
implying, in view of $\theta\in[0,1]$, the second claim of Theorem~\ref{ell1_risk} is proved.
%Note that (\ref{final})  reads
%$$
%\|f_1-f_2\|_1\geq O(1)T^{-1/2}\theta^{-1}\left[\sqrt{n\ln(2n/\beta)}+\sqrt{\ln(1/\epsilon)}\right].
%$$

\subsection{Proof of Proposition \ref{prop:corr} (Temporal Correlation of Sequential Detection Statistics)}

\begin{figure}[H]
\centering
\includegraphics[width = 0.7\textwidth]{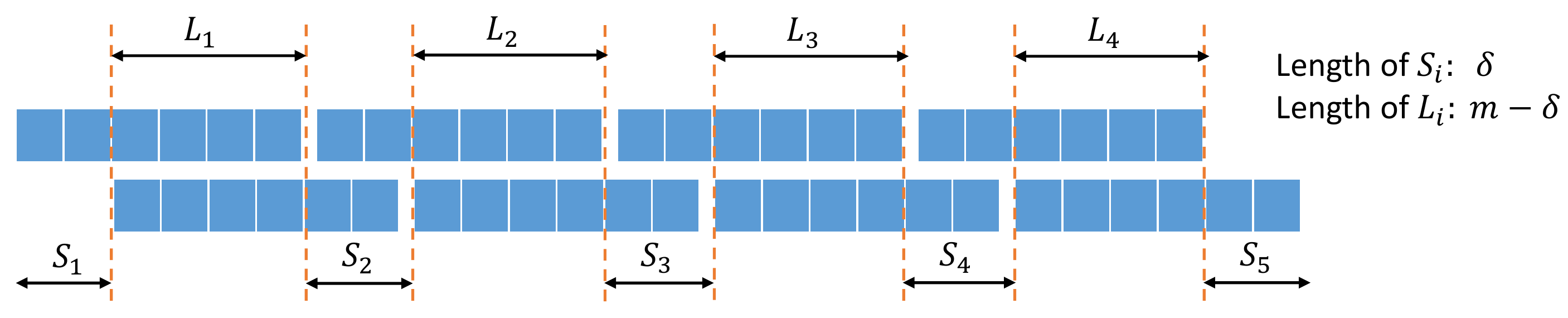} 
\caption{Sliding window illustration.}\label{fig:slidewindow}
\vspace{-0.1in}
\end{figure}

In order to approximate the correlation between statistics $\chi_{t,k}$ and $\chi_{\tau,s}$, we consider a simple case when the sample size $t-k=\tau-s=2m$. Denote the cardinality of the non-overlapping part in two time windows as $\delta = |t-\tau| = |k-s|$. Without loss of generality, we consider two sequences, $\{x_1,x_2,\ldots,x_{4m} \}$ and $\{x_{1+\delta},x_2,\ldots,x_{4m+\delta} \}$, with $4m-\delta$ overlapping elements, as illustrated in Fig.~\ref{fig:slidewindow}.
% {\it We focus on the setting that $\delta \ll m$, since the local dependency structure of the detection statistics plays an important role in characterizing the Average Run Length.}

Recall that $\mathrm{Var}(\chi_{t,k}) = 4\left[ \sum_{i=1}^n \sigma_i^2p_i^2(1-p_i)^2 + \sum_{i\neq j}\sigma_i\sigma_j p_i^2p_j^2\right]$ from the computation in Section~\ref{app:proof1}. Therefore, it remains to compute $\bE[\chi_{t,k}\chi_{\tau,s}]$. For simplicity, we compute $m^2\bE[\chi_t\chi_{t+\delta}]$ by writing it as a summation of several indicator functions as follows:
\[
\begin{aligned}
& m^2\bE[\chi_{t,k}\chi_{\tau,s}] \\
=& \bE\Bigg\{ \left[ \sum_{i=1}^n\sigma_i \left(\ind_i\{S_1\} +\ind_i\{L_1\}  - \ind_i\{S_3\}-\ind_i\{L_3\} \right) \left(\ind_i\{S_2\} +\ind_i\{L_2\}  - \ind_i\{S_4\}-\ind_i\{L_4\} \right)  \right] \\
& \quad \cdot \left[ \sum_{i=1}^n\sigma_i \left(\ind_i\{L_1\} +\ind_i\{S_2\}  - \ind_i\{L_3\}-\ind_i\{S_4\} \right) \left(\ind_i\{L_2\} +\ind_i\{S_3\}  - \ind_i\{L_4\}-\ind_i\{S_5\} \right)  \right] \Bigg\},
\end{aligned}
\]
where $\ind_i\{\mathcal S\} = \sum_{k\in\mathcal S} \ind\{x_k = i\}$ for a given set $\mathcal S$. Since $x_i$ are i.i.d. random variables, we have $\bE(\ind_i\{\mathcal S\}) = |\mathcal S|p_i$,  $\mathrm{Var}(\ind_i\{\mathcal S\}) = |\mathcal S|p_i(1-p_i)$, and $\mathbb E(\ind_i\{\mathcal S\}\ind_j\{\mathcal S\}) = |\mathcal S|(|\mathcal S|-1)p_ip_j$, with $|\mathcal S|$ denotes the cardinality. More specifically, for the decomposition shown in Fig.~\ref{fig:slidewindow}, we have $|S_i| = \delta$ and $|L_i| = m-\delta$. Substitute these into the above formulation, we have
\[
\begin{array}{c}
m^2\bE[\chi_{t,k}\chi_{\tau,s}] = \left[ 4(m-\delta)^2-2\delta^2\right] \left[ \sum_{i=1}^n \sigma_i^2p_i^2(1-p_i)^2 + \sum_{i\neq j}\sigma_i\sigma_j p_i^2p_j^2\right]. 
\end{array}
\]
% We decompose the above formula into two parts: I and II; I contains $n$ terms with the $i$-th term being
% \[
% \begin{aligned}
% & \sigma_i^2 \mathbb E\Bigg\{\left(\ind_i\{S_1\} +\ind_i\{L_1\}  - \ind_i\{S_3\}-\ind_i\{L_3\} \right) \left(\ind_i\{S_2\} +\ind_i\{L_2\}  - \ind_i\{S_4\}-\ind_i\{L_4\} \right) \\
% &\quad \cdot \left(\ind_i\{L_1\} +\ind_i\{S_2\}  - \ind_i\{L_3\}-\ind_i\{S_4\} \right) \left(\ind_i\{L_2\} +\ind_i\{S_3\}  - \ind_i\{L_4\}-\ind_i\{S_4\} \right)  \Bigg\}
% \end{aligned}
% \]
% Note that $L_i$'s and $S_i$'s are independent, after substituting the mean and variance of Bernoulli random variables, the above term can be simplified as:
% \[
% \left[ 4(m-\delta)^2-2\delta^2\right] \sigma_i^2p_i^2(1-p_i)^2.
% \]
% Similarly, we have that the part $B$ contains $n(n-1)$ terms with the $ij$-th term being
% \[
% \begin{aligned}
% &\sigma_i\sigma_j\mathbb E\Bigg\{\left(\ind_i\{S_1\} +\ind_i\{L_1\}  - \ind_i\{S_3\}-\ind_i\{L_3\} \right) \left(\ind_i\{S_2\} +\ind_i\{L_2\}  - \ind_i\{S_4\}-\ind_i\{L_4\} \right) \\
% &\quad \cdot \left(\ind_j\{L_1\} +\ind_j\{S_2\}  - \ind_j\{L_3\}-\ind_j\{S_4\} \right) \left(\ind_j\{L_2\} +\ind_j\{S_3\}  - \ind_j\{L_4\}-\ind_j\{S_4\} \right)  \Bigg\}
% \end{aligned}
% \]
% By an important observation $\mathbb E(\ind_i\{\mathcal S\}\ind_j\{\mathcal S\}) = |\mathcal S|(|\mathcal S|-1)p_ip_j$, we have the above term equals to 
% \[
% \left[ 4(m-\delta)^2-2\delta^2\right]\sigma_i\sigma_j p_i^2p_j^2.
% \]
Therefore, the correlation of the statistic $\chi_{t,k}$ and $\chi_{\tau,s}$ is
\[
\mathrm{Corr}(\chi_{t,k}, \chi_{\tau,s}) = \frac{4(m-\delta)^2-2\delta^2}{4m^2} = 1 - 2\frac{\delta}{m} + \frac{\delta^2}{2m^2}.
\]
Substitute $\delta= |t-\tau| = |k-s|$ into the above equation then we complete the proof.

\subsection{Proof of Theorem~\ref{thm:ARL} (ARL Approximation)}
The proof is based on a general method for computing first passage probabilities first introduced in \cite{yakir1998new} and further developed in \cite{SiegmundYakir2000} and \cite{SiegmundYakirZhang2010}, and commonly used in similar problems \cite{xie2013sequential,li2015m,cao2018multi}. 
First of all, it is worth mentioning that the probability measure in the following proof always stands for the nominal case where all samples are from the same distribution $p$. We define the variable
\[
Z_\tau = \tau(\xi_\tau - \eta_\tau)^\intercal \Sigma (\xi'_\tau - \eta'_\tau),
\]
where $\xi_\tau, \eta_\tau, \xi'_\tau,\eta'_\tau$ are empirical distributions of four independent (non-overlapping) sequences with equal length $\tau$. Recall that 
$\bE[Z_\tau] = 0$ and $\mathrm{Var}(Z_\tau) = \sigma_p^2$. We denote the moment generating function as
\[
\psi_\tau(\theta) = \log\bE\left[ \exp\{\theta Z_\tau \} \right],
\]
and select $\theta=\theta_\tau$ by solving the equation $\dot{\psi}_\tau(\theta)=b$. Since $Z_\tau$ is defined by a function of $4\tau$ independent random samples, $\phi_\tau$ converges to a limit as $\tau\rightarrow\infty$ and $\theta_\tau$ converges to a limiting value, denoted by $\theta$. The transformed distribution for all sequences at a fixed center position $k$ and the window size $\tau$ is denoted by $\bP_k^\tau$ and is defined by
\[
d\bP_k^\tau = \exp\{\theta Z_k^\tau - \psi_\tau (\theta_\tau)\} d\bP,
\]
where $Z_k^\tau = \tau(\xi_{k,\tau} - \eta_{k,\tau})^\intercal \Sigma (\xi'_{k,\tau} - \eta'_{k,\tau})$ is the statistic for location $k$ and window size $\tau$, as indicated in Fig.~\ref{fig:slidingwindow}. 

Let 
\[
\ell(k,\tau) := \log(d\bP_k^\tau/d\bP) = \theta Z_k^\tau - \psi_\tau (\theta_\tau).
\]
Denote $D=\{(k,\tau): 0\leq k\leq m, \lceil m_0/2 \rceil\leq \tau \leq \lceil m_1/2 \rceil\}$ be the set of all possible windows in the scan. Let $A=\{\max_{(k,\tau)\in D} Z_k^\tau \geq b\}$ be the event of interests (the event $\{\cT' \leq m\}$), i.e., the procedures stop before time $m$. 

By measure transformation, we have
\begin{equation}\label{eq:change-measure}
\begin{aligned}
\bP (A) & = \sum_{(k,\tau)\in D} \bE \Bigg[ \exp[\ell(k,\tau)] \Big( \sum_{(k',\tau')\in D}\exp[\ell(k',\tau')]\Big)^{-1};A\Bigg] \\
& = \sum_{(k,\tau)\in D} \bE_k^\tau \Bigg[ \Big( \sum_{(k',\tau')\in D}\exp[\ell(k',\tau')]\Big)^{-1};A\Bigg] \\
& = \sum_{(k,\tau)\in D}e^{\tilde{\ell}(k,\tau) - \ell(k,\tau)} \times \bE_k^\tau \left[ \frac{\max_{k',\tau'}e^{\ell(k',\tau')-\ell(k,\tau)}}{\sum_{k',\tau'}e^{\ell(k',\tau')-\ell(k,\tau)} }e^{-\tilde{\ell}(k,\tau) -[ \max_{k',\tau'}\ell(k',\tau') - \ell(k,\tau)]};A\right]\\
& = e^{-\theta_\tau\dot{\psi}_\tau(\theta_\tau) + \psi_\tau(\theta_\tau)} \times \sum_{(k,\tau)\in D}\bE_k^\tau \left[ \frac{M(k,\tau)}{S(k,\tau)} e^{-\tilde{\ell}(k,\tau) -\log M(k,\tau)};A\right]
\end{aligned}
\end{equation}
where 
\[
\begin{array}{rcl}
\tilde{\ell}(k,\tau) & =& \theta_\tau [Z_k^\tau - \dot{\psi}_\tau(\theta_\tau)], \\
M(k,\tau) & =& \max_{k',\tau'} \exp\{\theta_\tau (Z_{k'}^{\tau'} - Z_k^\tau) \}, \\
S(k,\tau) & =&  \sum_{k',\tau'} \exp\{\theta_\tau (Z_{k'}^{\tau'} - Z_k^\tau) \}.
\end{array}
\]
Since $k,\tau$ are fixed in much of the following analysis, we suppress the dependence of the notation on $k,\tau$ and simply writet $\tilde{\ell},S,M$. Under certain verifiable assumptions \cite{SiegmundYakirZhang2010}, a localization lemma allows us to simplify the expectation 
\[
\bE_k^\tau\left[ \frac{M}{S} e^{-\tilde{\ell} -\log M};\tilde{\ell} + \log M \geq 0 \right]
\]
into a simpler form
\[
\frac{1}{\sqrt{2\pi \sigma_{\tau}^2}} \bE\left[ \frac{M}{S} \right],
\]
where $\sigma_{\tau}^2$ stands for the variance of $\tilde{\ell}$ under measure $\bP_k^\tau$. The reduction relies on the fact that for large $m$, the local processes $M$ and $S$ are approximately independent of the global process $\tilde{\ell}$. Such independence allows the above decomposition into the expectation of $M/S$ times the expectation involving $\tilde{\ell} + \log M$, treating $\log M$ essentially as a constant. 

We first consider the process $M$ and $S$ and derive the expectation $\bE[M/S]$ following \cite{SiegmundYakir2000}. The difference between $Z_{k'}^{\tau'}$ and $Z_k^\tau$ can be writted in the form
\[
\begin{aligned}
Z_{k'}^{\tau'} - Z_k^\tau & = \tau'(\xi_{k',\tau'} - \eta_{k',\tau'})^\intercal \Sigma (\xi'_{k',\tau'} - \eta'_{k',\tau'}) -  \tau(\xi_{k,\tau} - \eta_{k,\tau})^\intercal \Sigma (\xi'_{k,\tau} - \eta'_{k,\tau}) \\
& = \tau'\left[ (\xi_{k',\tau'} - \eta_{k',\tau'})^\intercal \Sigma (\xi'_{k',\tau'} - \eta'_{k',\tau'}) - (\xi_{k,\tau} - \eta_{k,\tau})^\intercal \Sigma (\xi'_{k,\tau} - \eta'_{k,\tau}) \right] \\
& \quad + (\tau'-\tau) (\xi_{k,\tau} - \eta_{k,\tau})^\intercal \Sigma (\xi'_{k,\tau} - \eta'_{k,\tau}) 
\end{aligned}
\]

Observe that one may let $\tau' = \tau$ and substitute $\theta = \lim_{\tau\rightarrow\infty}\theta_\tau$ for $\theta_\tau$ in the definition of the increments and still maintain the required level of accuracy. When $\tau'=\tau$, the second term in the above expression vanish and the first term consists of two terms that are highly correlated. As characterized in Proposition~\ref{prop:corr}, when $\tau'=\tau$, the covariance between the two terms is given by
\[
\begin{aligned}
\mathbf{Cov}(\theta_\tau Z_{k'}^{\tau'},\theta_\tau Z_k^\tau) & = \theta_\tau^2\bE[Z_{k'}^{\tau'} , Z_k^\tau] = \theta_\tau^2\sigma_p^2\left(1 - 2\frac{|k'-k|}{\tau} + \frac{|k'-k|^2}{2\tau^2}\right).
\end{aligned}
\]
When $\tau$ is large, we have that the correlation depends on the difference $|k'-k|$ in a linear form, which shows that we have the random walk in the change time $k$, and the variance of the increment equals to $2\theta_\tau^2\sigma_p^2/\tau$. Following \cite{SiegmundYakir2000}, we have
\[
\bE[M/S] = [\theta_\tau^2\sigma_p^2/\tau\nu([2\theta_\tau\sigma_p^2/\tau]^{1/2})]^2.
\]

Moreover, the process $\tilde{\ell}$ is zero-mean and has variance $\sigma_\tau^2 = \mathrm{Var}_k^\tau(\tilde{\ell})=\theta_\tau^2\ddot{\psi}(\theta_\tau)$ under the measure $\bP_k^\tau$. Substituting the result for the expectations in \eqref{eq:change-measure} yields 
\[
\bP(\cT' \leq m) = 2\sum_{\tau=\lceil m_0/2 \rceil}^{\lceil m_1/2 \rceil}(m-2\tau) e^{-\theta_\tau\dot{\psi}_\tau(\theta_\tau) + \psi_\tau(\theta_\tau)}[2\pi \theta_\tau^2\ddot{\psi}_\tau(\theta_\tau)]^{-1/2}[\theta_\tau^2\sigma_p^2/\tau\nu([2\theta_\tau\sigma_p^2/\tau]^{1/2})]^2.
\]
In the limiting case, $Z_k^\tau$ can be well approximated using Gaussian distribution $\mathcal N(0,\sigma_p^2)$. The moment generating function then becomes $\psi(\theta) = \theta^2\sigma_p^2/2$, and the limiting $\theta=b/\sigma_p^2$, as the solution to $\dot{\psi}(\theta) = b$. Furthermore, the summation term can be approximated by an integral, to obtain
\begin{equation}\label{eq:exceed-prob}
\begin{aligned}
\bP(\cT' \leq m) & = 2\sum_{\tau=\lceil m_0/2 \rceil}^{\lceil m_1/2 \rceil}(m-2\tau) e^{-b^2/(2\sigma_p^2) }[2\pi b^2/\sigma_p^2]^{-1/2}[b^2/(\tau\sigma_p^2)\nu([2b^2/(\tau\sigma_p^2)]^{1/2})]^2 \\
& \approx 4e^{-b^2/(2\sigma_p^2) }[2\pi b^2/\sigma_p^2]^{-1/2}[b^2/\sigma_p^2]^2\int_{m_0/m}^{m_1/m} \nu^2([4b^2/(mt\sigma_p^2) ]^{1/2})(1-t)dt/t^2.
\end{aligned}
\end{equation}
Here it is assumed that $m$ is large, but small enough that the right-hand side of \eqref{eq:exceed-prob} converges to $0$ when $b\rightarrow \infty$. Changing variables in the integrand, we can rewrite this approximation as
\begin{equation}\label{eq:exponential}
\bP\{\cT'\leq m\} \approx m \times 2e^{-b^2/(2\sigma_p^2) }[2\pi b^2/\sigma_p^2]^{-1/2}  [b^2/\sigma_p^2] \int_{[4b^2/(m_1\sigma_p^2)]^{1/2}}^{[4b^2/(m_0\sigma_p^2)]^{1/2}} y \nu^2(y)dy. 
\end{equation}
From the arguments in \cite{siegmund1995using, siegmund2008detecting}, we know that $\cT'$ is asymptotically exponentially distributed and is uniformly integrable. Hence if $\lambda$ denotes the factor multiplying $m$ on the right-hand side of \eqref{eq:exponential}, then for large $m$, in the range where $m\lambda$ is bounded away from $0$ and $+\infty$, $\bP\{\cT'\leq m\} - [1-\exp(-\lambda m)] \rightarrow 0$. Consequently, $\bE[\cT'] \approx 1/\lambda$, thereby we complete the proof. Here we omit some technical details needed to make the derivation rigorous. Those details have been described and proved in \cite{SiegmundYakirZhang2010}.

\subsection{Proof of Theorem~\ref{thm:EDD} (EDD Approximation)}

Recall that $\chi_{t,0}$ is defined in \eqref{eq:detect_stat_online}. For any time $t$, we have
\[
\bE_0[\chi_{t,0}] = \frac{t}{2}(p-q)^\intercal \Sigma (p-q),
\]
which grows linearly with respect to time.
%Denote 
%\[
%\Delta^2= (p-q)^\intercal \Sigma (p-q)
%\] Let 
%\[
%\tilde{S}_t := \sum_{i=1}^t z_i,
%\]
%be a random walk where the increments $z_i$ are independent and identically distributed with mean $\Delta^2/2$ and variance $\Delta^2$. Define a stopping time $\tau = \min\{t: \tilde{S}_t > 0\}$. Our approximation to the expected detection delay depends on two related quantities. The first is
%\[
%\rho(\Delta) = \frac12 \bE[\tilde{S}_t^2]/\bE[\tilde{S}_t]
%\] 
%for which exact computational expressions and useful approximations are available in \cite{Siegmund1985}. In particular, we have
%\[
%\rho(\Delta) = \bE[z_1^2] / (2\bE[z_1])  - \sum_{i=1}^\infty i^{-1} \bE[\tilde{S}_i^{-}] = \Delta^2/4 + 1 - \sum_{i=1}^\infty i^{-1} \bE[\tilde{S}_i^{-}]
%\]
%where $(x)^{-} := -\min\{0,x\}$. The second quantity is $\bE[\min_{t\geq 0}\tilde{S}_t]$ which according to \cite{Siegmund1985}, is given by
%\[
%\bE[\min_{t\geq 0}\tilde{S}_t] = \rho(\Delta) - 1 - \Delta^2/4.
%\]
At the stopping time $\cT = T$, the expectation of the window-limited statistic in \eqref{eq:stoptime_window} can be computed if $m_1$ is sufficiently large (at least larger than the expected detection delay):
\[
\bE[\max_{0\leq k\leq T} \chi_{t,k} ] \approx \bE[ \chi_{t,0} ] =  \frac{T}{2}(p-q)^\intercal \Sigma (p-q).
\] On the other hand, we have that
\[
\bE[ \max_{0\leq k\leq T} \chi_{t,k} ]  = b + \bE[ \max_{0\leq k\leq T} \chi_{t,k} -b].
\] 
If we ignore the overshoot of the threshold over $b$ since it is a of order $o(b)$ when $b\rightarrow \infty$ (detailed analysis for overshoot has been developed in \cite{Siegmund1985}), then we obtain a first-order approximation as $b\rightarrow \infty$, by solving
\[
\frac{\bE_0[\cT]}{2}(p-q)^\intercal \Sigma (p-q) = b(1+o(1)).
\]
Therefore, a first-order approximation for the expected detection delay is given by
\[
\bE_0[\cT] = \frac{b(1+o(1))}{(p-q)^\intercal \Sigma (p-q)/2}.
\]
%The technical details and more rigorous statements have been developed in \cite{siegmund1995using,xie2013sequential,pollak1975approximations}. 

\subsection{Proof of Proposition~\ref{prop:dual} (Dual Reformulation)} 
For the minimization problem that defines $f(\sigma):=\min_{S\in{\cal S}}\Tr(\Sigma S)$, we introduce matrix variable $T\in\bR^{n\times n}$ such that $T_{ij} \geq |S_{ij}|, 1\leq i, j \leq n$ and $\sum_{i=1}^n \sum_{j=1}^n T_{ij} \leq 4$. We will jointly minimize over $S$ and $T$ in $f(\sigma)$.
The Lagrangian function of the problem $f(\sigma)$ writes
\[
\begin{array}{rcl}
&&\mathcal L(\sigma,S,P,T,\lambda,W,U,\xi, \{x_k \}_{k=1}^K)  \\
&= & \Tr(\Sigma S) - \Tr(PS) - r\Tr(SJ) + \lambda(\rho^2- \Tr(S)) \\
   && + \sum_k x_k(\Tr(SQ_k) - 4) + \xi(\Tr(TJ) - 4) - \Tr(U(T-S)) - \Tr(W(T+S)) \\
&= & \Tr\left( \left(\Sigma-P-rJ -\lambda I_n + \sum_{k=1}^Kx_k Q_k + U - W \right) S\right)  \\
&& + \lambda \rho^2 - 4\sum_{k=1}^K x_k - 4\xi + \Tr \left((\xi J - U - W)T\right),
\end{array}
\]
where $J \in \mathbb{R}^{n \times n}$ is a matrix with all elements equal to $1$, and we have
\[
f(\sigma) = \min_{S,T}\max_{\substack{ P \succcurlyeq 0, U \geq 0, W \geq 0 \\ \lambda \geq0, \xi \geq 0, x_k\geq 0,  r\in\bR}}\mathcal L(\sigma,S,P,T,\lambda,W,U,\xi, \{x_k \}_{k=1}^K). 
\]
%Use the fact that for any given $(y,t)$,
%\[
%\max_{(u,\lambda): \left\Vert u \right\Vert_2 \leq \lambda} - u^Ty-t\lambda = \max_{\lambda \geq 0} \lambda (\left\Vert y\right\Vert_2 - t) =
%\begin{cases}
%0, & \text{if $\left\Vert y\right\Vert_2 \leq t$} \\
%+\infty, & \text{otherwise}
%\end{cases}
%\]
%we can simplify the formulation in \eqref{eq:f2Lagrangian} as follows
%\[
%f^2(\sigma) = \min_{z\geq 0}\max_{\substack{ \omega_i \geq 0, \xi \geq 0, \lambda \geq0, \left\Vert s_k\right\Vert \leq x_k \\  \begin{bmatrix} p_i & q_i \\ q_i & r_i \end{bmatrix} \succcurlyeq 0 }} \mathcal L(\sigma,S,P,T,\lambda,W,U,\xi, \{x_k \}_{k=1}^K).
%\]
Then we have the dual problem of $f(\sigma)$ can be represented as
\[
\begin{array}{rcl}
\textrm{max}\quad && \lambda \rho^2 - 4\sum_k x_k  - 4\xi \\
\textrm{s.t.} \quad&& \lambda \geq0, P \succcurlyeq0, \xi \geq0, x_k \geq0, U \geq 0, W \geq 0,r\in\bR,\\
&& \sum_k x_kQ_k +U-W-P-rJ-\lambda I_n \succcurlyeq -\Sigma, \\
&& U_{ij} + W_{ij} \leq \xi, 1\leq i\leq n, 1\leq j\leq n. 
\end{array}
\]
Next we derive the dual form for $g(\sigma)$, similarly, we write the Lagrangian function as
\[
\begin{array}{c}
\mathcal L(\Sigma,P,\Lambda,V,\nu,\{\mu_k\}_{k=1}^K)  = \Tr((\Sigma^2 + \Lambda +V)P) + \nu(\Tr(PJ) - 1)+ \sum_k \mu_k(1-\Tr(PQ_k)),
\end{array}
\]
and we have
\[
g(\sigma) = \max_{P} \min_{\Lambda \succcurlyeq 0, V \geq 0, \nu \in \bR, \mu_k \geq 0} \mathcal L(\Sigma,P,\Lambda,V,\nu,\{\mu_k\}_{k=1}^K).
\]
Then we have the dual form of $g(\sigma)$ is
\[
\begin{array}{rcl}
\textrm{min}\quad && \sum_k \mu_k - \nu\\
\textrm{s.t.} \quad&& \Lambda \succcurlyeq 0, V \geq 0, \mu_k \geq 0, 1\leq k \leq K, \\
&& - \Lambda-V + \sum_k \mu_k Q_k - \nu J \succcurlyeq \Sigma^2.	
\end{array}
\]
Then the constraint $g(\sigma) \leq 1$ can be simplified to: $\exists \Lambda\succcurlyeq0, V\geq0, \mu_k\geq0,  \nu\in\bR$, such that
\[
\begin{array}{c}
\sum_k \mu_k - \nu \leq 1, \, - \Lambda -V + \sum_k \mu_k Q_k - \nu J \succcurlyeq \Sigma^2.
\end{array}
\]
Combine with the dual form of $f(\sigma)$, we have the problem \eqref{eq:optimal_weights} is equivalent to
\[
\begin{array}{rcl}
\textrm{max}\quad && \lambda \rho^2 - 4\sum_k x_k  - 4\xi \\
\textrm{s.t.} \quad&& \lambda \geq0, P \succcurlyeq0, \xi \geq0, x_k \geq0, U \geq 0, W \geq 0,r\in\bR,\\
&& \sum_k x_kQ_k +U-W-P-rJ-\lambda I_n \succcurlyeq -\Sigma, \\
&& U_{ij} + W_{ij} \leq \xi, 1\leq i\leq n, 1\leq j\leq n,\\
&& \sum_k \mu_k - \nu \leq 1, \\
&& \Lambda \succcurlyeq 0, V \geq 0, \mu_k \geq 0, 1\leq k \leq K, \\
&& - \Lambda -V + \sum_k \mu_k Q_k - \nu J \succcurlyeq \Sigma^2.
\end{array}
\]

\end{document}